\documentclass{article}
\usepackage{amsfonts,amsmath,amssymb,amsbsy}
\usepackage{graphicx}
\usepackage{epsfig}
\usepackage{theorem}
\newtheorem{theorem}{Theorem}
\newtheorem{lemma}{Lemma}

\begin{document}

\begin{titlepage}

\title{\textbf{Box dimension of unit-time map near nilpotent singularity of planar vector field}}
\author{Lana Horvat Dmitrovi\' c and Vesna \v Zupanovi\'c}

\maketitle

\begin{abstract}
The connection between discrete and continuous dynamical systems thro\-ugh the unit-time map has shown a significant role in bifurcation theory. Recently, it has also been used in fractal analysis of bifurcations. We study fractal properties of the unit-time map near nilpotent nonmonodromic singularities of planar vector fields using normal forms. We are interested in nilpotent singularities because they are nonhyperbolic, and we know that near nonhyperbolic singularities the box dimension is nontrivial. We study discrete orbits generated by the unit-time map, on the separatrices at the bifurcation point, and get results for the box dimensions of these orbits. Box dimension results will be illustrated in details using the examples of Bogdanov-Takens bifurcation and bifurcation of the nilpotent saddle. The study is also been extended to the appropriate singular points at infinity of normal form for the nilpotent singularity. Moreover, we study the unit-time map of the normal form for a saddle, near  singular points at infinity.
\end{abstract}


\end{titlepage} 

\def\b{\beta}
\def\g{\gamma}
\def\d{\delta}
\def\o{\omega}
\def\ty{\infty}
\def\e{\varepsilon}
\def\f{\varphi}
\def\:{{\penalty10000\hbox{\kern1mm\rm:\kern1mm}\penalty10000}}
\def\st{\subset}
\def\stq{\subseteq}
\def\q{\quad}
\def\M{{\cal M}}
\def\cal{\mathcal}
\def\eR{\mathbb{R}}
\def\eN{\mathbb{N}}
\def\Ze{\mathbb{Z}}
\def\Qu{\mathbb{Q}}
\def\Ce{\mathbb{C}}
\def\ov#1{\overline{#1}}
\def\D{\Delta}
\def\O{\Omega}

\def\bg{\begin}
\def\eq{equation}
\def\bgeq{\bg{\eq}}
\def\endeq{\end{\eq}}
\def\bgeqnn{\bg{eqnarray*}}
\def\endeqnn{\end{eqnarray*}}
\def\bgeqn{\bg{eqnarray}}
\def\endeqn{\end{eqnarray}}

\newcount\remarkbroj \remarkbroj=0
\def\remark{\advance\remarkbroj by1 \smallskip{Remark\ \the\remarkbroj.}\enspace\ignorespaces\,}

\pagestyle{myheadings}{\markright{Nilpotent singularities}}

\section{Introduction}

Since the 1970s fractal analysis has become a tool for studying dynamical systems, mostly their invariant sets and related measures. Complexity of invariant sets, measures, and graphs of  functions have been studied using numerous fractal dimensions. We stress widely used Hausdorff dimension and box dimension (also known as box counting, Minkowski dimension, Minkowski-Bouligand dimension, capacity dimension, limit capacity). This approach generated various results concerning box and Hausdorff dimensions of strange attractors, like Lorentz or Henon, and also results about Smale horseshoe, Julia and Mandelbrot sets. Results about homoclinic bifurcations and  fractal dimensions have been obtained, also. See survey article \cite{zuzu4} and references therein, to find mentioned results.

The new approach of studying dynamical systems by using the fractal dimension showed up recently. The idea is quite simple, 
to compute box dimension of any trajectory,  and to connect the obtained result with some other properties of the studied system. Results about box dimension and Minkowski content can be found in  recent development of  application of fractal analysis to  solutions of differential equations and dynamical systems. See for example \cite{lapo},  \cite{li}, \cite{pa1}, \cite{pazuzu}, \cite{zuzu3}. For our study of  discrete systems it is  particularly interesting fractal analysis of bifurcations of discrete dynamical systems (see \cite{neveda},  \cite{laho}, \cite{laho2}), \cite{zuzu}, \cite{belg}, \cite{laho3}, also applied to continuous systems.
 Direct connection between the box dimension of trajectories of dynamical systems and the bifurcation of that system has been proved in the articles.

 The first article in that direction was inspired by book of C.\ Tricot \cite{t}.  A new insight to fractal dimensions
 could be found  in the book.  We learnt that nonrectifiability (infinite length) of the curve could be measured by box dimension near the point of accumulation.  There are two interesting formulas in that book, which we exploited very much.  Formula for box dimension of nonrectifiable spiral 
 $r={\varphi}^{-\alpha}$, $0<\alpha \le1$, and formula for box dimension of nonrectifiable chirp $ f(x)=x^\alpha \sin x^{-\beta}$, $0<\alpha\le\beta$.
 Tricot's result about spiral authors applied to spiral trajectories of planar vector field in  article \cite{zuzu}.  
  Box dimension of spiral trajectories near weak focus or limit cycle is related  to the Hopf and Hopf-Takens bifurcation.  
  In that article it is observed that box dimension is related to the "potential" of the system to produce limit cycles under a small perturbation.  This observation made some kind of connection between box dimension and famous 16th Hilbert problem, so we think that this subject deserves to be explored. The study was extended involving Poincar\' e map (first return map) near focus or limit cycle.  The Poincar\' e map generates one-dimensional discrete system, so results  about  Poincar\' e map  from \cite{belg} are based on results about  discrete systems from \cite{neveda}.
  
Here we are interested in discrete systems, so we stress some results about discrete systems.  
First of all, notice that, for an orbit of one-dimensional discrete dynamical system, the Hausdorff dimension fails to show difference between systems. Namely, because of its property of countable stability, the Hausdorff dimension does not 'see' the countable sets at all. On the other hand, the box dimension is only finitely stable so it can 'see' them clearly. That is the reason why the box dimension is prefered for studying orbits of discrete dynamical systems. It is known that the box dimension of an orbit near hyperbolic fixed point of the one-dimensional discrete dynamical systems is trivial (see \cite{neveda}).  The analogous result for the hyperbolic fixed point of the systems in $\mathbb{R}^{n}$ is showed (see \cite{laho3}).  On the other hand, the box dimension near the nonhyperbolic fixed point of discrete dynamical system is strictly positive.  The article \cite{laho} proved the connection between the value of that box dimension and appropriate one and two-parameter bifurcations. Also in the article \cite{laho2} the box dimension result for the Neimark-Sacker bifurcation is proved. Recent work \cite{r} shows interesting connection beetween box dimension and Minkowski content with formal classification of parabolic diffeomorphisms.   
  
Discrete systems and continuous systems with spiral trajectories have been studied using box dimension.  
  Poincar\' e map makes a connection between continuous and discrete systems in the case  with spiral trajectories, we call it monodromic case.  For nonmonodromic case the connection could be established using unit-time map (also called time 1 map).
  Here we concentrate on nilpotent nonmonodromic singular points, while monodromic nilpotent singularities has been studied in   \cite{rezz}.
  In \cite{rzz} Hopf-Takens bifurcation at infinity has been studied using some generalization of box dimension and Poincar\' e compactification.  Expoiting this idea we also investigate unit-time map of nilpotent nonmonodromic singularities near singular points at infinity.
  More about dynamics near nilpotent singularities and near infinity could be found in  \cite{aga}, \cite{calito}, \cite{gi}, \cite {sz}, \cite{wu}, etc.

Now we recall the notions of box dimension and Minkowski content.  For further details see e.g.  \cite{fa}, and for some generalizations see  \cite{mrz}, \cite{rzz}.
%\cite{da1}.
Let $A\st\eR^N$ be bounded.  {\it Minkowski sausage} of radius $\e$ around $A$ is a $\e$-neighborhood of $A$, that is
$A_\e=\{y\in\eR^N\:d(y,A)<\e\}$. \ Let $s\ge0$.
\textit{The lower and upper $s$-dimensional Minkowski contents of $A$} are defined by
$$\M_*^s(A):=\liminf_{\e\to0}\frac{|A_\e|}{\e^{N-s}},\,\,\,\,\M^{*s}(A):=\limsup_{\e\to0}\frac{|A_\e|}{\e^{N-s}}.$$
Then \textit{the lower and upper box dimension} are defined by 
$$\underline\dim_BA=\inf\{s>0:\M_*^s(A)=0\}, \,\,\,\,\ov\dim_BA=\inf\{s>0:\M^{*s}(A)=0\}.$$
If $\underline\dim_BA=\ov\dim_BA$ we denote it by $\dim_BA$.
If there exists $d\ge0$ such that\ $0<\M_*^d(A)\le \M^{*d}(A)<\ty,$ then we say that set $A$ is \textit{Minkowski nondegenerate}.  Clearly, then $d=\dim_B A$.  If $|A_\e|\simeq \e^{s}$ for $\e$ small, then $A$ is Minkowski nondegenerate set and 
$\dim_B A=N-s$.  If $\M_*^s(A)=\M^{*s}(A)=\M^d(A)\in(0,\ty)$ for some $d\ge0$, then 
$A$ is said to be {\it Minkowski measurable}.  Clearly, then $d=\dim_BA$. \
Let $A$ and $B$ be two disjoint bounded sets such that $\dim_{B}A=\dim_{B}B$.  It is easy to see, using the finite stability of the upper box dimension
($\ov\dim_B(A\cup B)=\max\{\ov\dim_B A,\ov\dim_B B\}$), the monotonicity of the lower box dimension ($\underline\dim_B (A\cup B)\geq \underline\dim_B A$) and $\underline\dim_B (A\cup B)\leq \ov\dim_B (A\cup B)$ (for details see \cite{fa}), that 
\begin{equation} \label{jed1}
\dim_{B} (A\cup B) = \dim_{B} A =\dim_{B} B.
\end{equation}
In the paper the following definitions are used.
We say that any two sequences $(a_n)_{n\ge1}$ and $(b_n)_{n\ge1}$ of positive real numbers are {\it comparable} and write $a_n\simeq b_n$ as $n\to\ty$ if $A\le a_n/b_n\le B$ for some $A,B>0$  and  $n$ sufficiently big.
Analogously, two positive functions $f,g:(0,r)\rightarrow \eR$ are comparable and we write $f(x)\simeq g(x)$ as 
$x \rightarrow 0$ if $f(x)/g(x)\in [A,B]$ for $x$ small enough.

Hence we consider a discrete dynamical system 
$$\mathbf{x}_{n+1}=\mathbf{F}(\mathbf{x}_{n}),\, \mathbf{x}_1\in \mathbb{R}^{N}$$
generated by a $C^{k}$ function $\mathbf{F}:\eR^{N}\mapsto \eR^{N}$.  The \textit{orbit} of a system is a sequence $(\mathbf{x}_{n})_{n\geq1}$ such that $\mathbf{x}_{n+1}=\mathbf{F}(\mathbf{x}_{n})$ for some $\mathbf{x}_1\in \mathbb{R}^{N}$.   
Let $\mathbf{x}_0=0$ be a fixed point ($\mathbf{F}(\mathbf{x}_0)=\mathbf{x}_0$) of that system and let $A$ be a Jacobi matrix $D\mathbf{F}(\mathbf{x}_0)$ at $\mathbf{x}_0$.  The eigenvalues $\lambda_1,\ldots,\lambda_{N}$ of the matrix $A$ are called the {\it{multiplicators}} of fixed point.
We denote by $N_{0}$ the number of multipliers on the unit circle, by $N_{-}$ the number of multipliers inside the unit circle
and by $N_{+}$ the number of multipliers which lies outside the unit circle.   The fixed point is
{\it{hyperbolic}} if $N_{0}=0$, that is, there is no multipliers on the unit circle.  Hyperbolic point is called a hyperbolic saddle if $N_{-} N_{+}\neq 0$.  The fixed point is {\it{nonhyperbolic}} if $N_{0}\neq 0$.

In this paper the main object of our study is a box dimension of the orbit of the unit-time map of planar vector fields near the nilpotent singularity.  The unit-time map generates a two-dimensional discrete dynamical system which correspond to appropriate planar vector field.  
We know that the unit-time map on the characteristic orbit near the hyperbolic singularity is zero, while near the nilpotent singularity is positive and we will see that it depend on the order of a system and the asymptotic expansion of separatrices.
 
The remainder of this paper is organized as follows.  In Section 2 we recall the possible topological types of nilpotent singularities in the plane, and introduce the model system which we will study.  Then we present the unit-time maps for a model system with the nilpotent singularity. In Section 3 we apply  quasihomogenous blow-up method to the model system, and  get the asymptotic behavior of the separatrices, in the cases of nilpotent singularities which have the separatrices.  In Section 4 we  present the main result about the box dimension of the orbits of the unit-time map on the separatrices of nilpotent singularities.  Section 5  gives some examples of the bifurcations in the nilpotent singularities such as cusp and saddle.  For the cases of nilpotent saddle-node, node and the nilpotent singularity with hyperbolic and elliptic sector, we give one example of the appropriate system in order to illustrate the results in the previous chapters.
In Section 6 we study singularities at infinity for nilpotent model system, while in Section 7 we give some additional remarks about singularities which are not nilpotent. With the same technique we study normal form for saddle. We introduce a new notion of dual box dimension, which is related to dual Lyapunov constants, also called saddle quantities, see \cite{joyal}. We show an example where the saddle quantities of saddle normal form at the origin are related to the unit-time map of singularities of the normal form at infinity.

\section{Unit-time map}

First we will recall the known theorem about the classification of nilpotent singularities for planar vector fields (see \cite{dla}). We omitted  trivial cases (1) and (2).

\begin{theorem} {\rm\cite{dla}} \label{type}(\textbf{Nilpotent Singular Points})\\
Let $(0,0)$ be an isolated singular point of the vector field $X$ given by
\begin{eqnarray}
\dot{x}&=&y+A(x,y),\nonumber\\
\dot{y}&=&B(x,y),
\end{eqnarray}
where $A$ and $B$ are analytic functions in a neighborhood of the point $(0,0)$, and $j_{1}A(0,0)=j_{1}B(0,0)=0$.  Let $y=f(x)$ be the solution of the equation of  $y+A(x,y)=0$ in a neighborhood of the point $(0,0)$ and consider $F(x)=B(x,f(x))$ and $G(x)=(\frac{\partial A}{\partial x}+\frac{\partial B}{\partial y})(x,f(x))$.  Then the following holds:\\
(3) If $G(x)\equiv 0$ and $F(x)=ax^{m}+o(x^{m})$ for $m \in \eN$, $m\geq 1$, $a\neq 0$, then\\
(i) If $m$ is odd and $a>0$, then the origin is a saddle; and if $a<0$, then it is a center or a focus;\\
(ii) If $m$ is even then the origin is a cusp. \\
(4) If $F(x)=ax^{m}+o(x^{m})$ and $G(x)=bx^{n}+o(x^{n})$, $m,n\in\eN$, $m\geq 2$, $n\geq 1$, $a\neq 0$, $b\neq 0$, then we have\\
(i) If $m$ is even, and\\
(i1) $m<2n+1$, then the origin is a cusp. \\
(i2) $m>2n+1$, then the origin is a saddle-node. \\
(ii) If $m$ is odd and $a>0$, then the origin is a saddle. \\
(iii) If $m$ is odd, $a<0$ and\\
(iii1) Either $m<2n+1$, or $m=2n+1$ and $b^2+4a(n+1)<0$, then the origin is a center or a focus. \\
(iii2) $n$ is odd and either $m>2n+1$, or $m=2n+1$ and $b^2+4a(n+1)\geq 0$, then the phase portrait of the origin consist of one hyperbolic and one elliptic sector;\\
(iii3) $n$ is even and $m>2n+1$ or $m=2n+1$ and $b^2+4a(n+1)\geq 0$, then the origin is a node.  ($b>0$ repelling, $b<0$ attracting).
\end{theorem}

We are  interested in nonmonodromic isolated nilpotent singularities, it means that the singularities have separatrices or characteristic orbit going through the nilpotent singularity. The types are: cusp, saddle, saddle-node and node, and a nilpotent singularity with one elliptic and one hyperbolic sector.  

We consider the model system:
\begin{eqnarray} \label{sys2}
\dot{x}&=&y \nonumber\\
\dot{y}&=&f(x)+yg(x)+y^2B(x,y),
\end{eqnarray}
where $f$, $g$ and $B(x,y)$ are $C^{\infty}$ functions, $j_{1} f(0)=g(0)=j_{\infty}B(0,0)=0$.
Also it is satisfied $j_{\infty}f(0)\neq 0$, moreover $f(x)=ax^{m}+o(x^{m})$, $a\neq 0$.  For $g$ there are two cases: $j_{\infty}g(0)= 0$ or $g(x)=bx^{n}+o(x^{n})$, $b\neq 0$.
So we consider the system
\begin{eqnarray} \label{sys2}
\dot{x}&=&y\nonumber\\
\dot{y}&=&ax^{m}+bx^{n}y+y^2B(x,y)+o(x^{m})+y o(x^{n})
\end{eqnarray}
under the assumption  $\deg(B)+2>\max\{m,n+1\}$.

Before proving the theorems, we recall the procedure for calculating the unit-time map of continuous system by using the Picard iterations.
So, we consider the continuous dynamical system
\begin{equation} \label{k1}
\dot{\mathbf{x}}=\mathbf{F}(\mathbf{x})
\end{equation}
where $\mathbf{x}\in\mathbb{R}^{N}$, $F:\mathbb{R}^{N}\rightarrow \mathbb{R}^{N}$.  The simpliest way of getting the discrete dynamical system from the continuous one is by using the unit-time map
 $\phi_{t}(\mathbf{x})$.  Namely, we fix $t_0>0$ and we consider the system  which is generated by the iteration of the map $\phi_{t_0}$ (map with displacement $t_0$ along the trajectory of ($\ref{k1}$)).  If we take $t_0=1$, we get the discrete dynamical system generated by the unit-time map
\begin{equation} \label{k2}
\mathbf{x} \mapsto \phi_{1}(\mathbf{x}).
\end{equation} 
It can be easily shown that the isolated fixed points of ($\ref{k2}$) corresponds to the isolated singularities of ($\ref{k1}$).  
In order to study the connection between the hyperbolicity and stability of these points, we need to find the connection between the corresponding eigenvalues of $D\mathbf{F}(\mathbf{x}_0)$ and $D\phi_{1}(\mathbf{x}_0)$.

Now we look at the continuous dynamical system with the singularity $\mathbf{x_0}=0$ 
\begin{equation} \label{k3}
\mathbf{\dot{x}}=\mathbf{F}(\mathbf{x})=A\mathbf{x} + \mathbf{F}^{(2)}(\mathbf{x})+\mathbf{F}^{(3)}(\mathbf{x})+\ldots, \,\,\,\mathbf{x}\in\mathbb{R}^{N},
\end{equation}
where $A=D\mathbf{F}(0)$ and $\mathbf{F}^{(k)}$ are smooth polinomial vector function of order  $k$: $\mathbf{F}^{(k)}(\mathbf{x})=\mathcal{O}(\left\|\mathbf{x}\right\|^{k})$ i $$F_{i}^{(k)}(\mathbf{x})=\sum_{j_1+\ldots+j_{n}=k} b_{i,j_1,\ldots,j_{n}}x_1^{j_1}
x_2^{j_2}\ldots x_{n}^{j_{n}}.$$
We denote the corresponding flow of ($\ref{k1}$) with $\phi_{t}(\mathbf{x})$.  Now we would like to find Taylor expansion of $\phi_{t}(\mathbf{x})$  near $\mathbf{x_0}=0$ by using the process of Picard iterations.  Namely, let $\mathbf{x}^{(1)}(t)=e^{At}\mathbf{x}$ be a solution of linear equation  $\dot{\mathbf{x}}=A\mathbf{x}$ with the initial value  $\mathbf{x}$, and define  
$$\mathbf{x}^{(k+1)}(t)=e^{At}\mathbf{x} + \int_0^{t} e^{A(t-\tau)}(\mathbf{F}^{(2)}(\mathbf{x}^{(k)}(\tau))+\ldots + \mathbf{F}^{(k+1)}(\mathbf{x}^{(k)}(\tau)))d\tau.$$
It is easy to show that $(k+1)$-iteration does not change the terms of order $l\leq k$.  By the substitution $t=1$ in $\mathbf{x}^{(k)}(t)$ we get the Taylor expansion of the unit-time map $\phi_1(\mathbf{x})$ until the terms of order $k$ 
\begin{equation}
\phi_1(\mathbf{x})=e^{A}\mathbf{x} + \mathbf{g}^{(2)}(\mathbf{x})+\ldots + \mathbf{g}^{(k)}(\mathbf{x}) + \mathcal{O}(\left\|\mathbf{x}\right\|^{k+1}),
\end{equation}
where $\mathbf{g}^{(i)}$ are polynomial vector function of a form as functions $\mathbf{F}^{(i)}$.
We get $B=D\phi_1(0)=e^{A}$, where $A=D\mathbf{F}(0)$.  It means that $\mathbf{x_0}=0$ is a hyperbolic (nonhyperbolic) singularity of $(\ref{k1})$ if and only if  $\mathbf{x}_0=0$ is a hyperbolic (nonhyperbolic) fixed point of map ($\ref{k2}$).  In dimension one, it is obvious because $e^0=1$.  In the plane, we have three cases: 
\begin{itemize}
\item $A$ has two different real eigenvalues $\lambda_1\neq\lambda_2$ $\Rightarrow$ $B$ has two different real eigenvalues $e^{\lambda_1}$ and $e^{\lambda_2}$
\item $A$ has one real eigenvalue $\lambda$ $\Rightarrow$ $B$ has one real eigenvalue $e^{\lambda}$
\item $A$ has two complex conjugated eigenvalues $\lambda_{1,2}=a\pm bi$ $\Rightarrow$ $B$ has two complex conjugated eigenvalues $e^{a}(\cos b \pm i\sin b )$ 
\end{itemize}

Using the above procedure, we can prove the following lemma.

\begin{lemma} ({\textbf {The unit-time map}})\label{flow}\\
Let $(x_0,y_0)=(0,0)$ be a nilpotent singularity of the  system $(\ref{sys2})$.  Then the following holds:
\begin{enumerate}
	\item If $m<n+1$ then the unit-time map has a form
		\begin{eqnarray} 
\hskip-0.7cm		x_{k+1}&=& x_{k}+y_{k}+\frac{a}{2} x_k^{m}+ac_{11}x_k^{m-1}y_k+\ldots+ac_{1m}y_k^{m}+O(\left\|x\right\|^{m+1})\nonumber\\
	\hskip-0.7cm	y_{k+1}&=&y_{k}+ax_k^{m}+ad_{11}x_k^{m-1}y_k+\ldots+ad_{1m}y_k^{m}+ O(\left\|x\right\|^{m+1});
		\end{eqnarray}
\hskip-0.3cm		with the constants $c_{1i}=c_{1i}(m)$, $d_{1i}=d_{1i}(m)$.
	\item If $m=n+1$ then the unit-time map has a form
	\begin{eqnarray}
\hskip-0.3cm	x_{k+1}&=&x_{k}+y_{k}+\frac{a}{2}x_k^{m}+ac_{21}x^{m-1}y_k+\ldots+ac_{2,m-1}x_ky_k^{m-1}+O(\left\|x\right\|^{m+1})\nonumber\\
\hskip-0.3cm	y_{k+1}&=&y_{k}+ax_k^{m}+ad_{21}x_k^{m-1}y_k+\ldots+ad_{2m}y_k^{m}+ O(\left\|x\right\|^{m+1});
	\end{eqnarray}
\hskip-0.3cm		with the constants $c_{2i}=c_{2i}(m)$, $d_{2i}=d_{2i}(m)$.
	\item If $m>n+1$ then the unit-time map has a form
	\begin{eqnarray} \hskip-0.3cm x_{k+1}\hskip-0.2cm&=&\hskip-0.2cmx_{k}+y_{k}+\frac{b}{2}x_{k}^{n}y+bc_{31}x_{k}^{n-1}y_{k}^{2}+\ldots+bc_{3,n-1}x_{k}y_{k}^{n}+\frac{b}{n+2}y_{k}^{n+1}+O(\left\|x\right\|^{n+2})\nonumber\\
\hskip-0.2cm	y_{k+1}\hskip-0.2cm&=&\hskip-0.2cmy_{k}+bx_{k}^{n}y_{k}+bd_{31}x_{k}^{n-1}y_{k}^2+\ldots+bd_{3,n-1}x_{k}y_{k}^{n}+\frac{b}{n+1}y_{k}^{n+1}+ O(\left\|x\right\|^{n+2});
	\end{eqnarray}
	\hskip-0.5cm	with the constants $c_{3i}=c_{3i}(n)$, $d_{3i}=d_{3i}(n)$, $i=1,\ldots,n-1$.
\end{enumerate}
\end{lemma}

\textbf{Proof.}\\
In the case $m\leq n+1$, by using the above procedure we find the Taylor expansion of the unit-time map up to terms of order $m$, while in the case
$m>n+1$ we can get the Taylor expansion up to  terms of order $n+1$.$\blacksquare$\\

\section{Separatrices}

We continue with the model system:
\begin{eqnarray} \label{sys3}
\dot{x}&=&y \nonumber\\
\dot{y}&=&f(x)+yg(x)+y^2B(x,y)
\end{eqnarray}
where $f$, $g$ and $B(x,y)$ are $C^{\infty}$ functions, $j_{1} f(0)=g(0)=j_{\infty}B(0,0)=0$.
We study the systems with $j_{\infty}f(0)\neq 0$, that is $f(x)=ax^{m}+{\rm{o}}(x^{m})$, $a\neq 0$.  For $g$ there are two cases: $j_{\infty}g(0)= 0$ or $g(x)=bx^{n}+{\rm{o}}(x^{n})$, $b\neq 0$.

In order to find the asymptotics of the separatrices we  use  quasihomogenous blow-up (see \cite{dla}, \cite{z}). In the quasihomogenous blowing-up procedure we make a finite sequence of changes of variables in the system, leading to the desingularisation of the system. We choose new variables with respect to the Newton diagram of the system, it depends on the slope of the side in the diagram. For nilpotent singularities it is possible to desingularise   vector field using only one quasihomogenous blow-up. Regarding Newton diagram, we have 3 different cases for nilpotent singularities:\\
\begin{enumerate}
	\item Hamiltonian like case ($m<2n+1$)
	\item Singular like case ($m>2n+1$)
	\item Mixed case $m=2n+1$
\end{enumerate}

Now we would like to find the asymptotics of the separatrices for each case.\\

\textbf{1. Case $m<2n+1$}\\
The system ($\ref{sys3}$) is $C^{\infty}$ is equivalent to
\begin{eqnarray} \label{sys4}
\dot{x}&=&y \nonumber\\
\dot{y}&=&\delta x^{m}+y(bx^{n}+\rm{o}(x^{n}))+\rm{O}(y^2).
\end{eqnarray}

\textbf{\underline{Case 1A} $m$ even:} \\
For $\delta=1$ we  have a cusp.
Using the homogeneous blow-up:
\begin{eqnarray}
x&=&u^2\nonumber\\
y&=&u^{m+1}\bar{y}
\end{eqnarray}
we get the system
\begin{eqnarray}
\dot{u}&=&\frac{1}{2}u \bar{y}\nonumber\\
\dot{\bar{y}}&=&(1-\frac{m+1}{2}\bar{y}^2)+O(u),
\end{eqnarray}
with the singularities: $(u,\bar{y})=(0,\pm \sqrt{\frac{2}{m+1}})$, which are hyperbolic saddles.  
Now we find the Taylor approximations for the stable and unstable manifolds of this systems, in order to find the approximation for the separatrices of the initial system.  We look at the unstable manifold for the singularity $(0,\sqrt{\frac{2}{m+1}})$, and the stable one for $(0,-\sqrt{\frac{2}{m+1}})$.  The invariant unstable manifold is
$$\bar{y}=\sqrt{\frac{2}{m+1}}+\alpha_1 u +\alpha_2 u^2+O(u^3).$$
The stable manifold is
$$\bar{y}=-\sqrt{\frac{2}{m+1}}+\beta_1 u +\beta_2 u^2+O(u^3),$$
and the corresponding separatrices are:\\
unstable $$y=\sqrt{\frac{2}{m+1}}x^{\frac{m+1}{2}}+\alpha_1 x^{\frac{m+2}{2}} +\alpha_2 x^{\frac{m+3}{2}}+O(x^{\frac{m+4}{2}});$$
stable $$y=-\sqrt{\frac{2}{m+1}}x^{\frac{m+1}{2}}+\beta_1 x^{\frac{m+2}{2}} +\beta_2 x^{\frac{m+3}{2}}+O(x^{\frac{m+4}{2}}).$$
So, the asymptotic behavior of the separatrices near the origin (cusp) is $$y\simeq \pm \sqrt{\frac{2}{m+1}}x^{\frac{m+1}{2}}.$$
This is what we need in order to calculate the box dimension of the unit time map.  

\textbf{Remark 1.} Since the asymptotics depends only on $m$, it is interesting to explore where is hidden the power $n$.
In fact, it is easy to see that the second term in the asymptotic series of separatrix depend on $n$.  
If we include the separatrices in the system, we get that $\alpha_1=\ldots=\alpha_{k-1}=0$ and $\alpha_{k}\neq 0$ for $k=2n+1-m$.  For example, if $m=2$ and $n=1$, we get $y=\sqrt{\frac{2}{3}}x^{\frac{3}{2}}+\alpha_1 x^2+\alpha_2x^{\frac{5}{2}}+\ldots$, while for $m=2$ and $n=3$, we have 
$y=\sqrt{\frac{2}{3}}x^{\frac{3}{2}}+\alpha_5 x^4+\alpha_{10} x^{\frac{13}{2}}+\ldots$. \

\textbf{\underline{Case 1B} $m$ odd:}\\
For $\delta=1$, we have a nilpotent saddle.
The case $\delta=-1$, where we have a center or focus, we do not study here. 
Using the homogenuous blow-up:
\begin{eqnarray}
x&=&u\nonumber\\
y&=&u^{\frac{m+1}{2}}\bar{y}
\end{eqnarray}
we get the system
\begin{eqnarray}
\dot{u}&=&u \bar{y}\nonumber\\
\dot{\bar{y}}&=&(1-\frac{m+1}{2}\bar{y}^2)+O(u),
\end{eqnarray}
with the singularities: $(u,\bar{y})=(0,\pm \sqrt{\frac{2}{m+1}})$, which are hyperbolic saddles.  
Then, we find Taylor approximations for the stable and unstable manifolds of this systems, in order to find the approximation for the separatrices of initial system.  We look at the unstable manifold for the singularity $(0,\sqrt{\frac{2}{m+1}})$, and the stable one for $(0,-\sqrt{\frac{2}{m+1}})$.  The invariant unstable manifold is
$$\bar{y}=\sqrt{\frac{2}{m+1}}+\alpha_1 u +\alpha_2 u^2+O(u^3),$$
and the stable manifold 
$$\bar{y}=-\sqrt{\frac{2}{m+1}}+\beta_1 u +\beta_2 u^2+O(u^3).$$
Corresponding separatrices are:\\
unstable $$y=\sqrt{\frac{2}{m+1}}x^{\frac{m+1}{2}}+\alpha_1 x^{\frac{m+3}{2}} +\alpha_2 x^{\frac{m+5}{2}}+O(x^{\frac{m+7}{2}});$$
stable $$y=-\sqrt{\frac{2}{m+1}}x^{\frac{m+1}{2}}+\beta_1 x^{\frac{m+3}{2}} +\beta_2 x^{\frac{m+5}{2}}+O(x^{\frac{m+7}{2}}).$$
So, the asymptotic behaviiour of the separatrices near the nilpotent saddle is $$y\simeq \pm \sqrt{\frac{2}{m+1}}x^{\frac{m+1}{2}}.$$
Notice that it is the same as in the cusp case.

\textbf{Remark 2.} It is easy to see that the second term in the asymptotic series depend on $n$.  
If we include the separatrices in the system, we get that the separatices are:
$$y=\pm \sqrt{\frac{2}{m+1}}x^{\frac{m+1}{2}}+\alpha_{k}x^{\frac{m+1}{2}+k}+\alpha_{2k}x^{\frac{m+1}{2}+2k}+\ldots$$ where $k=\frac{2n+1-m}{2}$.
For example, for $m=3$ and $n=2$, we get $y=\sqrt{\frac{1}{2}}x^{2}+\alpha_1 x^3+\alpha_2x^{4}+\ldots$. \

\textbf{2. Case $m>2n+1$:}\\
The system ($\ref{sys3}$) is $C^{\infty}$ equivalent to 
\begin{eqnarray} \label{sys5}
\dot{x}&=&y \nonumber\\
\dot{y}&=&a x^{m}+y(x^{n}+o(x^{n}))+O(y^2), \q a\neq 0
\end{eqnarray}
Nilpotent singularities in this case are nilpotent saddle, saddle-node, node and singularity with the elliptic and hyerbolic sector.
Using the homogenuous blow-up:
\begin{eqnarray}
x&=&u\nonumber\\
y&=&u^{n+1}\bar{y}
\end{eqnarray}
we get the system
\begin{eqnarray}
\dot{u}&=&u \bar{y}\nonumber\\
\dot{\bar{y}}&=&\bar{y}(1-(n+1)\bar{y}+O(u))+au^{m-2n-1}
\end{eqnarray}
with two singularities: $T_1(u,\bar{y})=(0,\frac{1}{n+1})$ and $T_2(u,\bar{y})=(0,0)$.  The singularity $T_1$ is a hyperbolic saddle, while $T_2$ is a semi-hyperbolic point with the unstable manifold $u=0$, and the center manifold transverse to it.
The semi-hyperbolic singularity $T_2$ can be saddle, node or saddle-node. Topological type of the singularity of the system (\ref{sys5}) depends of the type of $T_2$.

 Now, we look at the unstable manifold for the singularity $T_1(0,\frac{1}{n+1})$.  The invariant unstable manifold is
$$\bar{y}=\frac{1}{n+1}+\alpha_1 u +\alpha_2 u^2+O(u^3),$$
so the corresponding unstable separatrix has a form
 $$y=\frac{1}{n+1}x^{n+1}+\alpha_1 x^{n+2}+\alpha_2 x^{n+3}+O(x^{n+4}).$$
Then, the asymptotic behavior of the "upper" separatrix near the origin is $$y\simeq \frac{1}{n+1}x^{n+1}.$$
This is what we need in order to calculate the box dimension of the unit-time map.  

\textbf{Remark 3.} Since the asymptotics depends only on $n$, it is interesting to explore where is hidden the power $m$.
In fact, it is easy to see that the second term in the asymptotic series depend on $m$.  
If we include the separatrices in the system, we get that $\alpha_1=\ldots=\alpha_{k-1}=0$ and $\alpha_{k}=\frac{n+1}{m-n}$, for $k=2n+1-m$, that is,
$$y=\frac{1}{n+1}x^{n+1}+\alpha_{m-(2n+1)} x^{m-n}+O(x^{m-n+1}).$$

Now, we  find the center manifold for the singularity $T_2=(0,0)$.  We know that has a form
$$\bar{y}=\gamma_1 u +\gamma_2 u^2+O(u^3),$$
so the corresponding invariant manifold for the initial system is\\
$$y=\gamma_1 x^{n+2}+\gamma_2 x^{n+3}+O(x^{n+4}).$$
But, if we include it in the system, we get that the first term  also depends on $m$.  Since $m>2n+1$, there exists $k>0$ such that $k=m-2n-1$.  We obtain the asymptotic series 
$$y=\gamma_{k}x^{n+k+1}+O(x^{n+k+2})=\gamma_{m-2n-1}x^{m-n}+O(x^{m-n+1}).$$

\textbf{3. Case $m=2n+1$:}\\
The system ($\ref{sys3}$) is $C^{\infty}$ is equivalent to 
\begin{eqnarray} \label{sys6}
\dot{x}&=&y \nonumber\\
\dot{y}&=&a x^{2n+1}+y(x^{n}+o(x^{n}))+O(y^2), \q a\neq 0.
\end{eqnarray}
In this case we are interested in the nilpotent singularities such as saddle, node and singularity with the elliptic and hyperbolic sector.  

Using the homogenuous blow-up:
\begin{eqnarray}
x&=&u\nonumber\\
y&=&u^{n+1}\bar{y}
\end{eqnarray}
we get the system
\begin{eqnarray}
\dot{u}&=&u \bar{y}\nonumber\\
\dot{\bar{y}}&=&a+\bar{y}-(n+1)\bar{y}^2+\bar{y}O(u)
\end{eqnarray}
with several possibilities:
\begin{enumerate}
\item For $a>0$, we have two singularities $T_{1,2}=(0,\frac{1\pm\sqrt{4a(n+1)+1}}{2(n+1)})$ which are both hyperbolic saddles, and the corresponding topological type of the initial system is a nilpotent saddle. \
We denote by $A_{1,2}=\frac{1\pm\sqrt{4a(n+1)+1}}{2(n+1)}$, where $A_1>0$ and $A_2<0$.  Then the unstable separatrix is of a form
$$y=A_1x^{n+1}+\alpha_1 x^{n+2}+\alpha_2 x^{n+3}+O(x^{n+4}),$$
and the stable is 
$$y=A_2x^{n+1}+\beta_1 x^{n+2}+\beta_2 x^{n+3}+O(x^{n+4}).$$

\item For $a<0$, we can have two singularities $T_{1,2}=(0,\frac{1\pm\sqrt{4a(n+1)+1}}{2(n+1)})$, where $T_1$ is a hyperbolic saddle and $T_2$ is a node; or if $4a(n+1)+1=0$), we have only one singularity $T_1=(0,\frac{1}{2(n+1)})$ which is a saddle-node.  By blowing down, we can get the nilpotent node (for $n$ even) or a nilpotent singularity with one elliptic and one hyperbolic sector (for $n$ odd). \
The invariant unstable separatrix for the singularity $T_1$ is:\\
if $1+4a(n+1)>0$, then
$$y=A_1x^{n+1}+\alpha_1 x^{n+2}+\alpha_2 x^{n+3}+O(x^{n+4});$$
if $1+4a(n+1)=0$, then
$$y=\frac{1}{2(n+1)}x^{n+1}.$$  

\item If $4a(n+1)+1<0$, then there is no singularities, and the initial system has a center or focus.
\end{enumerate}

We showed that in all nonmonodromic cases $y\simeq x^{n+1}$.

\section{Box dimension of the unit-time map}

As we mentioned before, it is already known that the box dimension of each orbit of discrete dynamical system near the hyperbolic fixed point in $\mathbb{R}^{N}$ is 0. Also, it is known that the box dimension of orbit near the nonhyperbolic fixed point of discrete system in $\mathbb{R}^{N}$ is strictly positive, see \cite{laho3}. These results can be applied to the unit-time map of the continuous system near nonhyperbolic singularity.

\begin{theorem} {\rm \cite{laho3}}
Let $(0,0)$ be a hyperbolic singular point of continuous planar dynamical system.  Then the unit-time map on each characteristic trajectory near $(0,0)$ has positive box dimension.
\end{theorem}

The results of box dimension in the case of the nonhyperbolic singularity with only one multiplier on the unit circle can be found in \cite{laho3}.
Now we would like to get the result for the exact value of box dimension near a nilpotent singularity.

\begin{lemma}
Let $A=\{x_{k}\}_{k\in\mathbb{N}}$ and $B=\{y_{k}\}_{k\in\mathbb{N}}$ be a two decreasing sequences which tends to $0$ with initial points $x_0$ and $y_0$ and with the properties $x_{k}-x_{k+1}\simeq x_k^{\alpha}$, 
for $\alpha>1$ and $y_{k}-y_{k+1}\simeq y_k^{\beta}$, for $\beta>1$. Let $S(x_0,y_0)=\{(x_{k},y_{k})\}$ be a two-dimensional discrete 
dynamical system, with initial point $(x_0,y_0)$.
Then the following holds:\\
(i) if $\alpha\geq\beta$, then $\dim_{B}S=1-\frac{1}{\alpha}$; \\
(ii) if $\alpha<\beta$, then $\dim_{B}S=1-\frac{1}{\beta}$. \\
\end{lemma}
\textbf{Proof.}\\
From Theorem 1, \cite{neveda}, it follows that $\dim_{B}A=1-\frac{1}{\alpha}$ and $\dim_{B}B=1-\frac{1}{\beta}$. It is obvious that the set $A$ is an orthogonal projection of the set $S$ on the $x$-axis. Analogously, the set $B$ is an orthogonal projection of $S$ on the $y$-axis. Since the orthogonal projection is a Lipscitz map, then we know that $\underline{\dim}_{B}S\geq \rm{max}\{\dim_{B}A,\dim_{B}B\}$. That is the lower bound for the box dimension of the set $S$. The upper bound for the box dimension can be calculated directly by estimating the area of the Minkowski sausage of radius $\varepsilon.$
We denote by $n_{A}(\varepsilon)$ the minimal $n\in\mathbb{N}$ for which $x_{n}-x_{n+1}<2\varepsilon$, and analogously for $n_{B}(\varepsilon)$. Also, we denote by $n_{S}(\varepsilon)$ the minimal $n\in\mathbb{N}$ such that $\sqrt{(x_{n}-x_{n+1})^2+(y_{n}-y_{n+1})^2}<2\varepsilon$.
It is easily seen that in the case $\alpha\geq\beta$, it holds
\begin{equation} \label{nejed1}
n_{A}(\sqrt{2}\varepsilon)\leq n_{S}(\varepsilon)\leq n_{B}(\sqrt{2}\varepsilon).
\end{equation}
In the case $\alpha<\beta$, the opposite inequalities are valid.
Now, we can estimate the area of the Minkowski sausage of radius $\varepsilon$ which we divide into the tail $\left|S_{\varepsilon}\right|_{t}$ (before overlapping), and nucleus $\left|S_{\varepsilon}\right|_{n}$ (after overlapping).
We want to calculate the upper bound for the box dimension. So we have
$$\left|S_{\varepsilon}\right|_{t}=\pi\varepsilon^2(n_{S}(\varepsilon)-1),$$
and
$$\left|S_{\varepsilon}\right|_{n}\leq \pi\varepsilon^2+\int_{0}^{x_{n_{S}(\varepsilon)}}(g(x)+\delta-(g(x)-\delta))dx$$
where $y=g(x)$ is the curve on which the discrete orbits lies, and it can be easily seen that $g(x)\simeq x^{\gamma}$, where $\gamma=\frac{\alpha-1}{\beta-1}$. Constant $\delta=\rm{max}_{n\in\mathbb{N}}\delta_{x_{n}}$ is a maximum $\delta$ (see Figure 1) such that almost the whole $S_{\varepsilon}$ (without two semicircles) is between the curves $g(x)+\delta$ and $g(x)-\delta$. In the case $\alpha\geq\beta$, the curve $y=g(x)$ is increasing, and the derivative is also increasing so $\delta=\delta_{x_{n_{S}(\varepsilon)}}$. \\
Now we have
$$\left|S_{\varepsilon}\right| \leq \left|S_{\varepsilon}\right|_{n}+\left|S_{\varepsilon}\right|_{t} = \pi\varepsilon^2 n_{S}(\varepsilon)+2\delta x_{n_{S}(\varepsilon)}.$$
Notice that $\delta_{x_{n}}\geq\varepsilon$ and that $\lim_{n \rightarrow \infty} \delta_{x_{n}}=\varepsilon$. 
So, we can choose $\varepsilon$ small enough such that $\delta<2\varepsilon$. Then we have
\begin{equation} \label{nejed2}
\left|S_{\varepsilon}\right| \leq \pi\varepsilon^2 n_{S}(\varepsilon)+4\varepsilon x_{n_{S}(\varepsilon)}.
\end{equation}
In the case $\alpha\geq\beta$, from $(\ref{nejed1})$ it follows 
$$x_{n_{S}(\varepsilon)}\leq x_{n_{A}(\sqrt{2}\varepsilon)},$$
that is,
$$ n_{S}(\varepsilon)^{-\frac{1}{\alpha-1}} \leq n_{A}(\sqrt{2}\varepsilon)^{-\frac{1}{\alpha-1}}.$$
Now we put $x_{n_{A}(\sqrt{2}\varepsilon)}\simeq n_{A}(\varepsilon)^{-\frac{1}{\alpha-1}}$ and $n_{A}(\varepsilon)\simeq \varepsilon^{-(1-\frac{1}{\alpha})}$ in the inequality $(\ref{nejed2})$, divide by $\varepsilon^{2-s}$, and get
$$\frac{\left|S_{\varepsilon}\right|}{\varepsilon^{2-s}} \leq \pi \varepsilon^{s} C_1 \varepsilon^{-(1-\frac{1}{\beta})} + 2C_2\varepsilon^{s-1} n_{A}(\varepsilon)^{-\frac{1}{\alpha-1}}\leq$$
$$\leq C_1\pi \varepsilon^{s-1+\frac{1}{\beta}} + 2C_2 \varepsilon^{s-1}\varepsilon^{-(1-\frac{1}{\alpha})(-\frac{1}{\alpha-1})}\leq$$
$$\leq C_1\pi\varepsilon^{s-(1-\frac{1}{\beta})}+C_2\varepsilon^{s-(1-\frac{1}{\alpha})}.$$
So it follows that $$\overline{\dim}_{B}S\leq 1-\frac{1}{\alpha},$$ and we proved the lemma.
$\blacksquare$\\ 

\textbf{Remark 4.} This lemma can also be easily proven by using the Lemma 4 from \cite{zuzu3q}. 

\textbf{Remark 5.} In general, $y=g(x)$ from the previous proof can be a piecewise spline which numerically aproximates the curve. 
Also, notice that the discrete orbits studied here lie on the curves, that is, trajectories or separatrices of the continuous system.

\begin{center}
\includegraphics[width=5cm,bb=0in 0in 2in 2in]{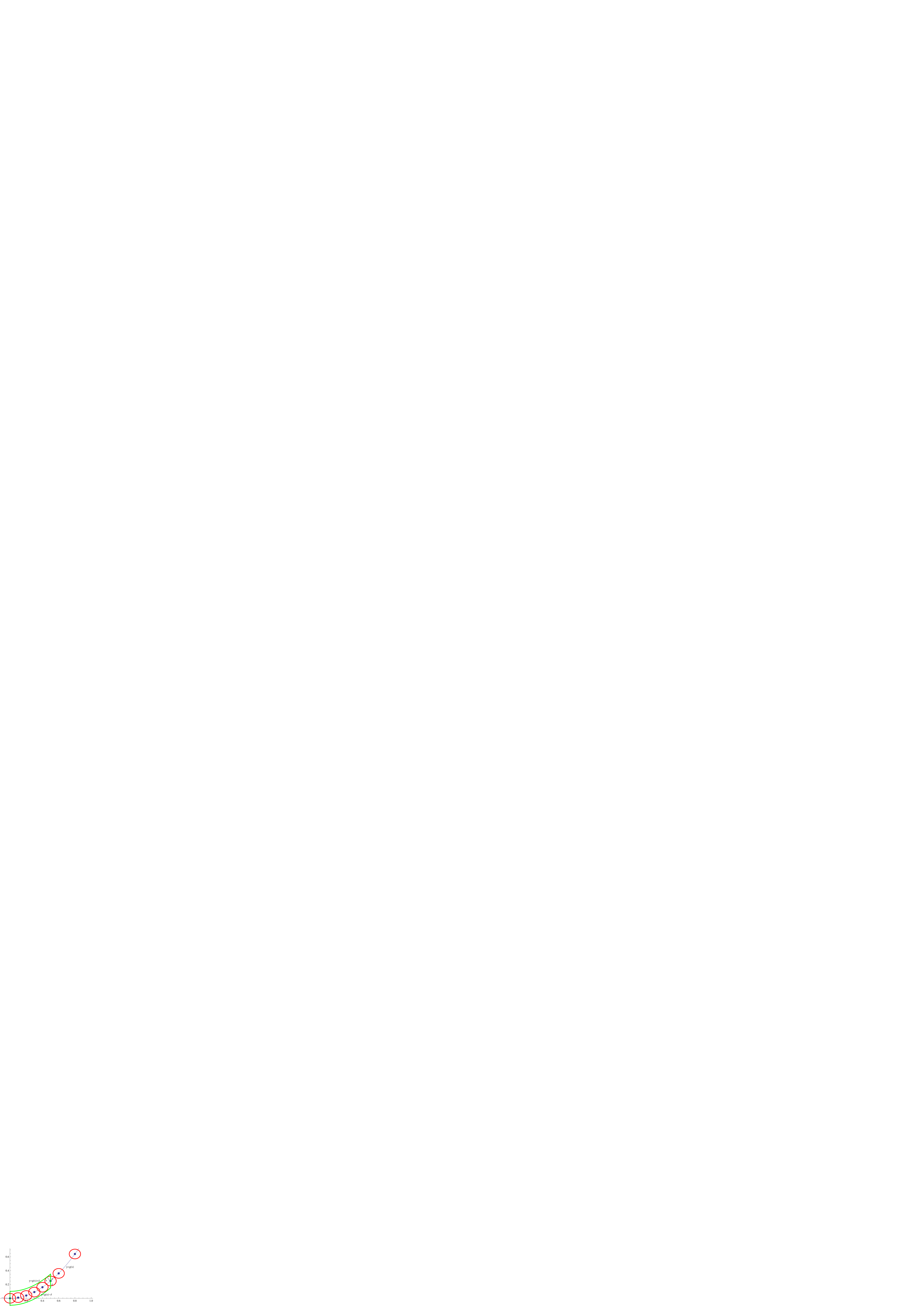}\\ 
\textbf{Figure 1} upper bound for $\left|S_{\varepsilon}\right|$, $\alpha\geq\beta$
\end{center}

\begin{theorem} \label{box}({\textbf{Box dimension near nilpotent singularity}})\\
Let we have a system (\ref{sys2}).  Let $(x_0,y_0)=(0,0)$ be a nilpotent singularity
and let $y(x)=x^{\gamma}+o(x^{\gamma})$, $\gamma \in (1,m)$ be a separatrix of the system (\ref{sys2}). 
Let $\Gamma$ be an orbit on the separatrix generated by the unit-time map of the system near $(0,0)$, and $\Gamma_x$, $\Gamma_y$ are projections of $\Gamma$ to the coordinate axes.
Then the following holds:\\
\textbf{ (1) } If $m\le n+1$, then $\dim_{B}\Gamma_{x}=1-\frac{1}{\gamma}$, $\dim_{B}\Gamma_{y}=1-\frac{\gamma}{m}$.\\
(i)If $\gamma^2\geq m$, then $\dim_{B}\Gamma=\dim_{B}\Gamma_{x}=1-\frac{1}{\gamma}$.\\ 
(ii) If $\gamma^2< m$, then $\dim_{B}\Gamma=\dim_{B}\Gamma_{y}=1-\frac{\gamma}{m}$. \\
\textbf{(2)} If $m>n+1$, then $\dim_{B}\Gamma_{x}=1-\frac{1}{\gamma}$, $\dim_{B}\Gamma_{y}=1-\frac{\gamma}{n+\gamma}$.\\ 
(i) If $\gamma\geq \frac{1}{2}(1+\sqrt{1+4n})$, then $\dim_{B}\Gamma=\dim_{B}\Gamma_{x}=1-\frac{1}{\gamma}$.\\ 
(ii) If $\gamma < \frac{1}{2}(1+\sqrt{1+4n})$, then $\dim_{B}\Gamma=\dim_{B}\Gamma_{y}=1-\frac{\gamma}{n+\gamma}$. 
\end{theorem}

\textbf{Proof.}\\
(1) Case $m\leq n+1$:\\
From Lemma 1, it follows that the asymptotics of the unit-time map are
$$x_{k}-x_{k+1}\simeq x_{k}^{\gamma},$$
and
$$y_{k}-y_{k+1}\simeq y_{k}^{\frac{m}{\gamma}}.$$
Using Theorem 1 from \cite{neveda}, we have 
$$\dim_{B}\Gamma_{x}=1-\frac{1}{\gamma},$$ 
and $$\dim_{B}\Gamma_{y}=1-\frac{\gamma}{m}.$$
We denote by $\alpha=\gamma$, $\beta=\frac{m}{\gamma}$, then we have
\begin{eqnarray}
\alpha\geq \beta \,\,\Leftrightarrow \,\,\gamma^2\geq m\nonumber\\
\alpha<\beta \,\,\Leftrightarrow \,\, \gamma^2< m.\nonumber
\end{eqnarray}
So the results for the box dimension follow from Lemma 2.\\
(2) Case $m>n+1$: \\
From Lemma 1, it follows that the asymptotics of the unit-time map are
$$x_{k}-x_{k+1}\simeq x_{k}^{\gamma},$$
and
$$y_{k}-y_{k+1}\simeq y_{k}^{\frac{n}{\gamma}+1}.$$
Using Theorem 1 from \cite{neveda}, we have 
$$\dim_{B}\Gamma_{x}=1-\frac{1}{\gamma}$$ and $$\dim_{B}\Gamma_{y}=1-\frac{\gamma}{n+\gamma}.$$
Denoting by $\alpha=\gamma$, $\beta=\frac{n}{\gamma}+1$, for $\gamma>0$ we get
\begin{eqnarray}
\alpha\geq \beta \,\,\Leftrightarrow \,\,\gamma\geq \frac{1}{2}(1+\sqrt{1+4n}),\nonumber\\
\alpha<\beta \,\,\Leftrightarrow \,\, \gamma< \frac{1}{2}(1+\sqrt{1+4n}).\nonumber
\end{eqnarray}
Now the results for the box dimension easily follow from Lemma 2.\\
$\blacksquare$\\

The results from the theorem will be illustrated by  examples in Section 5.

\textbf{Remark 6.} It is interesting to notice that for cusp and nilpotent saddle, we have a characteristic set of values concerning box dimension.
For cusp 
$$
D_c=\{{\frac{   \dim_{B}\Gamma_{x}  }{  \dim_{B}\Gamma_{y}   }}= {\frac{4k}{2k+1}}    :k\in\eN\}=\left\{\frac43,\,\frac85,\,\frac{12}7,\frac{16}9,\,\frac{20}{11},\,\dots\right\}.
$$
The set $D_c$ coincide to the set of values of box dimensions of the spiral trajectory near weak focus, see Theorem 9, \cite{zuzu}.
Also for the nilpotent saddle, the possible values are: 
$$
D_s=\{{\frac{   \dim_{B}\Gamma_{x}  }{  \dim_{B}\Gamma_{y}   }}=2-{\frac{1}{k+1}}:k\in\eN\}=\left\{\frac32,\,\frac83,\,\frac{7}4,\frac{9}5,\,\frac{11}{6},\,\dots\right\}.
$$
The set $D_s$ coincide to the set of values of box dimensions of the spiral trajectory near limit cycle, see Theorem 10, \cite{zuzu}.

\textbf{Remark 7.} In the nilpotent case the connection between the box dimension and the multiplicity of  fixed point, or the cyclicity of singularity should be further explored.  Cyclicity means maximal number of limit cycles which could be obtain from the system under small perturbation.
Unfoldings of nilpotent singularities show their complex structure after application of blowing-up method. Blowing-up the singularity we find polycycles "inside". Roughly speaking polycycles are separatrices which are "closed curves" passing through singularities, the simplest cases are saddle-loop and two saddle-loop. Poincar\' e map near singularity, limit cycle or polycycle is a standard tool for studying  cyclicity.
Problem is that Poincar\' e map near monodromic nilpotent singularity or polycycle is not analytic. These fractal methods could be adjusted to such cases using some other scale to obtain an asymptotic expansions, for Chebyshev scale see \cite{mrz}.

\section{Examples}

In this section we will present several examples of the nilpotent singularities in order to illustrate the results for box dimensions of the unit-time map on the separatrices.  In the case of cusp and nilpotent saddle, we will present the whole unfolding for the appropriate bifurcation, while in other cases we will give only the nilpotent situation.
In the following examples  we also use results from 
\cite{neveda} and \cite{laho} dealing with box dimension of the orbits of discrete one-dimensional systems at the bifurcation point. The corresponding values for nondegenerate saddle-node, and period doubling bifurcations are $\frac12$, $\frac23$, while for hyperbolic orbit of node and focus, box dimension is equal to $0$. For saddles box dimension is computed on the stable and unstable manifold. The semihyperbolic cases can be reduced to the center manifold, see \cite{laho3}.
To complete the study about box dimension of the whole unfolding we also use result about discrete Hopf bifurcation called Neimark-Sacker appearing in 2-dimensional systems, see \cite{laho2}.
Discrete spiral orbit at the bifurcation parameter has box dimension equal to $\frac43$.

\subsection{Cusp}

Let us consider the normal form for the Bogdanov-Takens bifurcation
\begin{eqnarray}
\dot{x}&=&y\nonumber\\
\dot{y}&=& \beta_1+\beta_2 x +x^2-xy,
\end{eqnarray}
where $\beta_{1,2}\in\mathbb{R}$ are parameters.  We can see the bifurcation diagram of Bogdanov-Takens bifurcation at Figure 2, see more details in \cite{kuz}.  
We denote by $H$ the negative part of $\beta_2$ axis, because it is a curve where Hopf bifurcation occurs.

For $\beta_1=\beta_2=0$ we have a cusp and we use Theorem \ref{box}, case (1), (i), with $m=2$, $n=1$, $\gamma=\frac{m+1}2=\frac32$, and get $\dim_{B}S=\dim_{B}S_{x}=1/3$, $\dim_{B}S_{y}=\frac{1}{4}$ on the separatrices. See Figure 3a. For other cases we use the results from \cite{neveda}, \cite{laho}, and \cite{laho2}. At region 1 there are no singularities. Furthermore, by passing through the curve $T-$ a saddle and a node appear, so it is a saddle-node bifurcation curve. On the center manifold we have $\dim_{B}S=\dim_{B}S_{x}=\frac{1}{2}$ and $\dim_{B}S_{y}=0$ (Figure 3c). Somewhere in the region 2 the node becomes a focus, and crossing the curve $H$ a limit cycle is born. So box dimension on the Hopf bifurcation curve is $\frac{4}{3}$ (Figure 3e). Passing through the curve $P$ saddle homoclinic bifurcation occurs, that is a saddle-loop is apeared (for recent result about $P$ see \cite{bogdanov}). In region 4 the saddle-loop is broken and there are two singularities, a saddle and a node. If we continue the journey clockwise and finally return to region 1, once more a saddle-node bifurcation occurs (curve $T^{+}$).

All such objects are unfolded in the cusp with  $\dim_{B}S=1/3$, for $\beta_1=\beta_2=0$. Notice that the box dimension is nontrivial when some local bifurcation occurs. To detect the global bifurcation on $P$, box dimension near homoclinic loop should be computed, see \cite{mrz}. All hyperbolic cases inside the regions have trivial box dimensions.
At the following figures, the trajectories of continuous system are drawn in blue colour, while the orbits of discrete system generated by the unit-time map are drawn in red. The separatrices are green. Accumulation of red points on the separatrix has been measured by the box dimension.

\begin{center}
\includegraphics[width=4cm,bb=0in 0in 2in 2in]{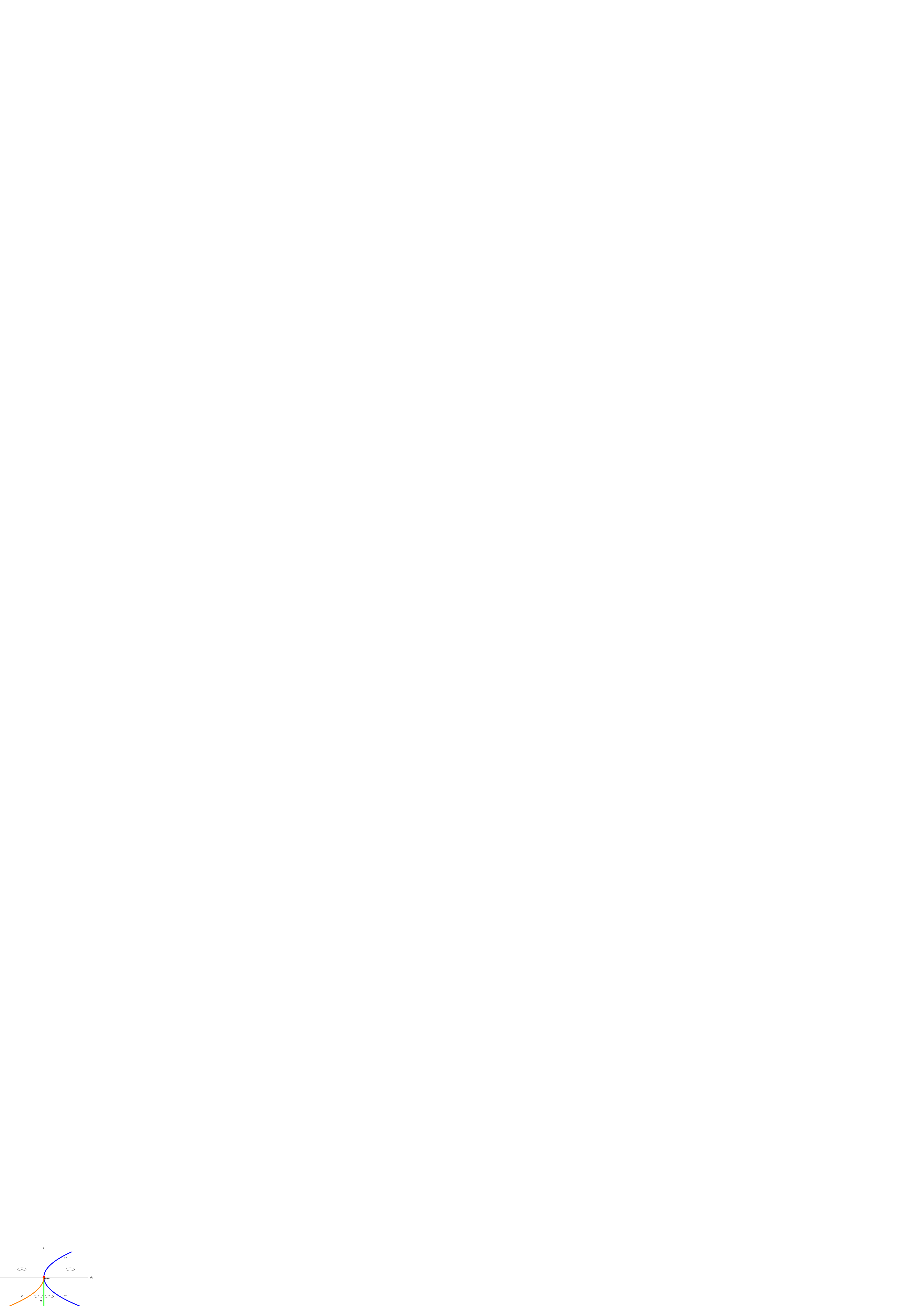} \\ 
 \textbf{Figure 2} bifurcation diagram 
\end{center}

\begin{center}
\includegraphics[width=3cm,bb=0in 0in 3in 3in]{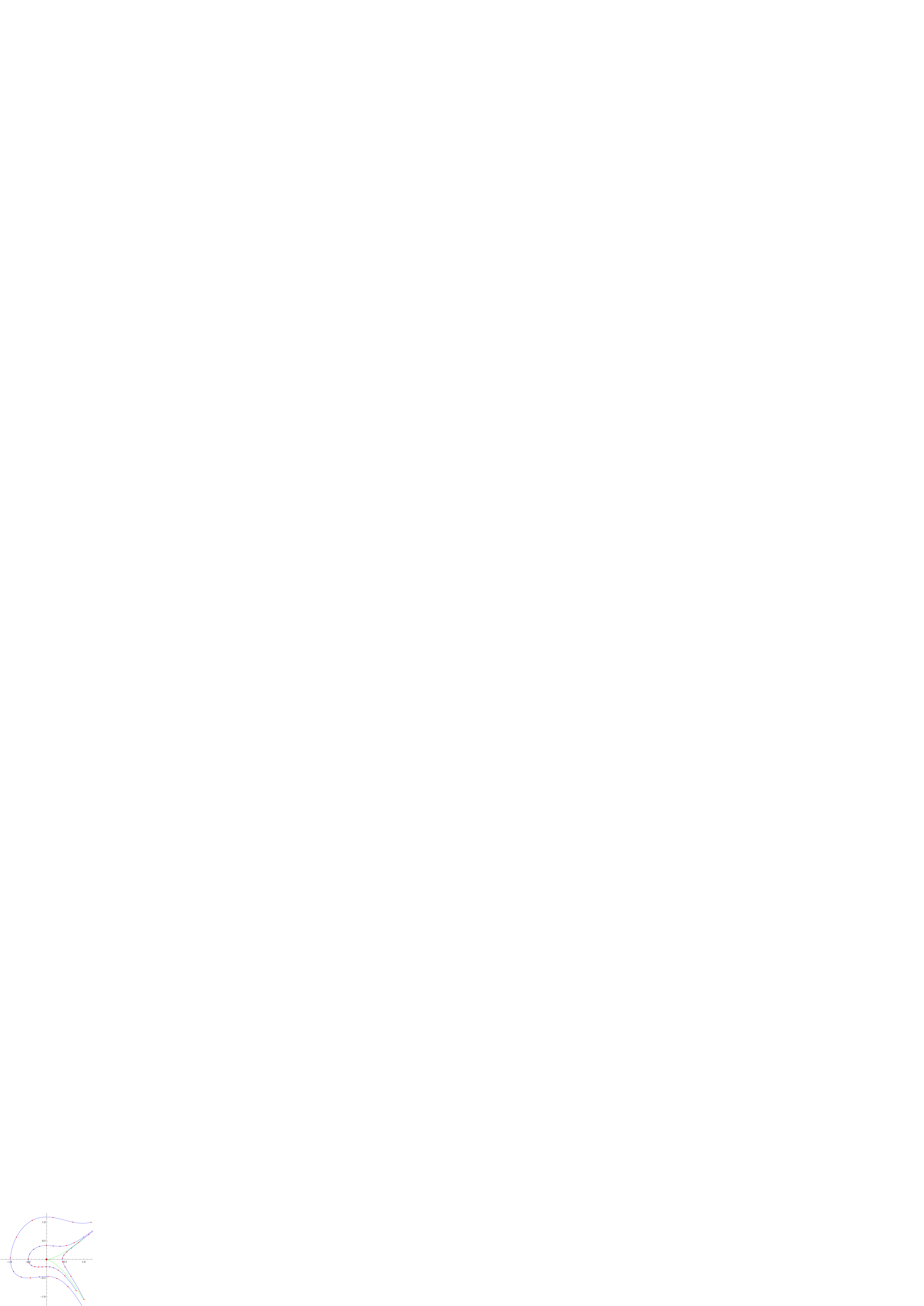} \hskip 2cm \includegraphics[height=3cm,bb=0in 0in 3in 3in]{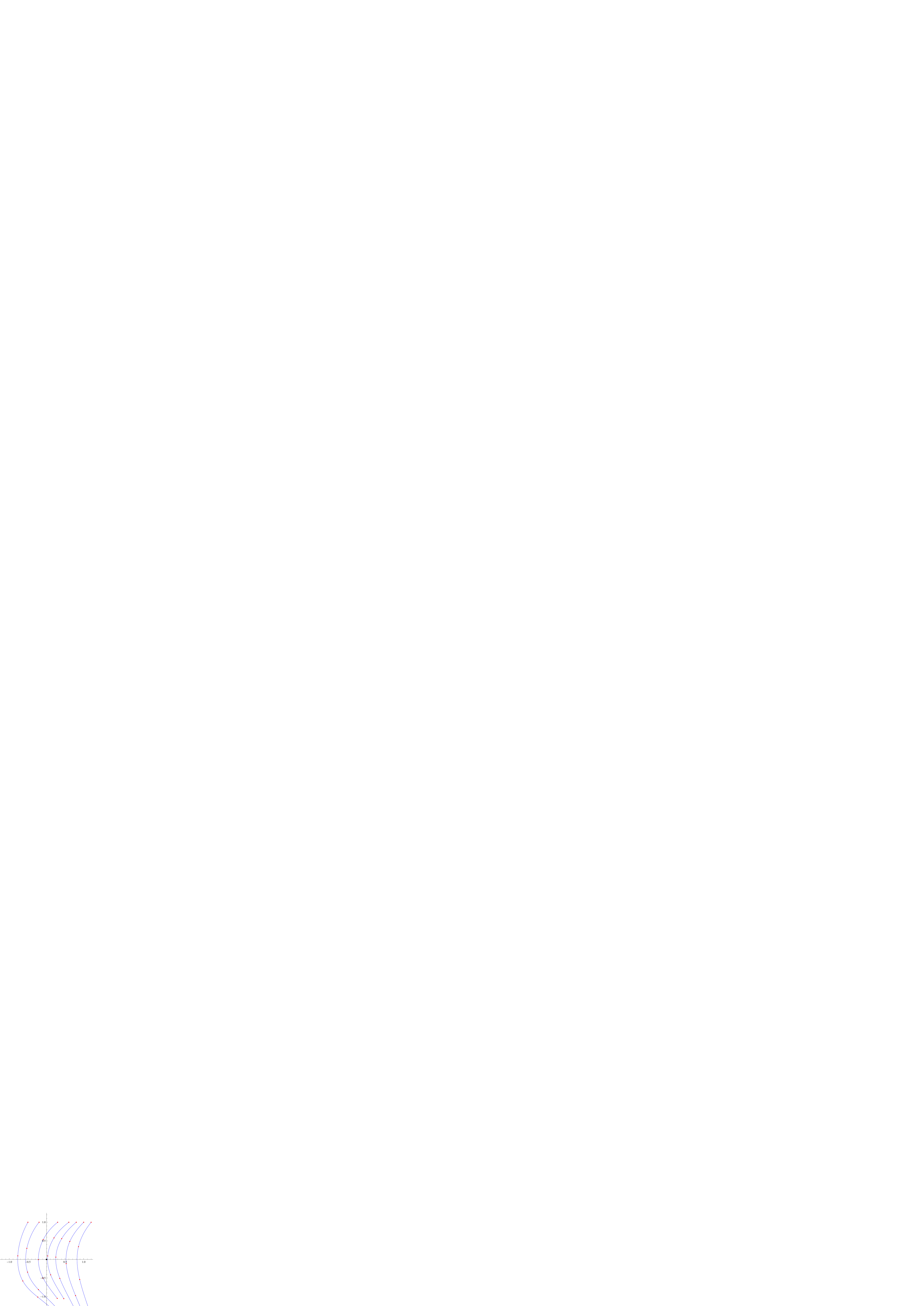}\\ 
\textbf{Figure 3a} cusp, $\beta_1=\beta_2=0,\dim_{B}S=1/3$\hskip 0.1cm \textbf{Figure 3b} region $1$ from Figure 2
\end{center}

\begin{center}
\includegraphics [width=4cm,bb=0in 0in 4.5in 4.5in]{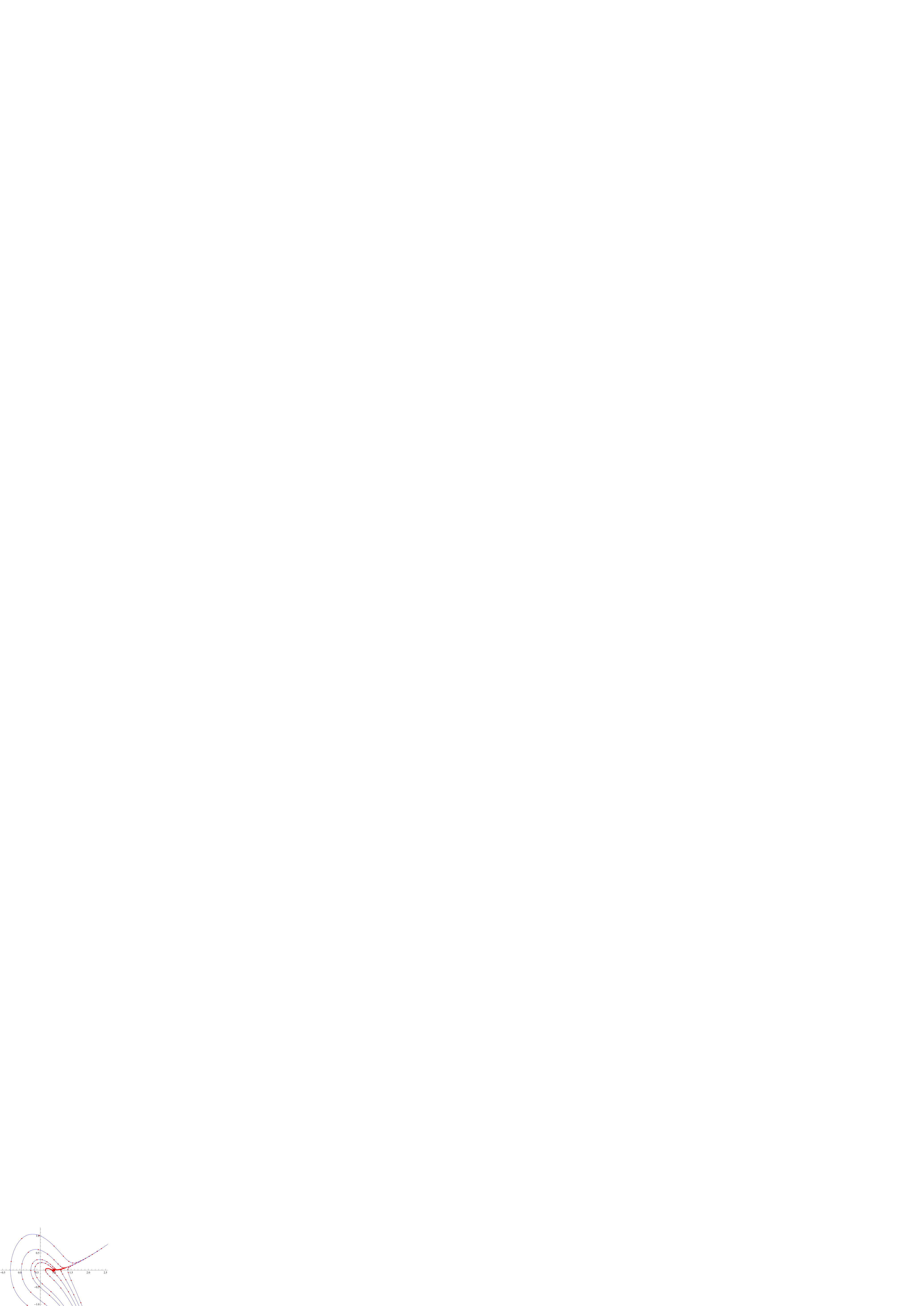} \hskip 2cm \includegraphics[height=3cm,bb=0in 0in 4.5in 4.5in]{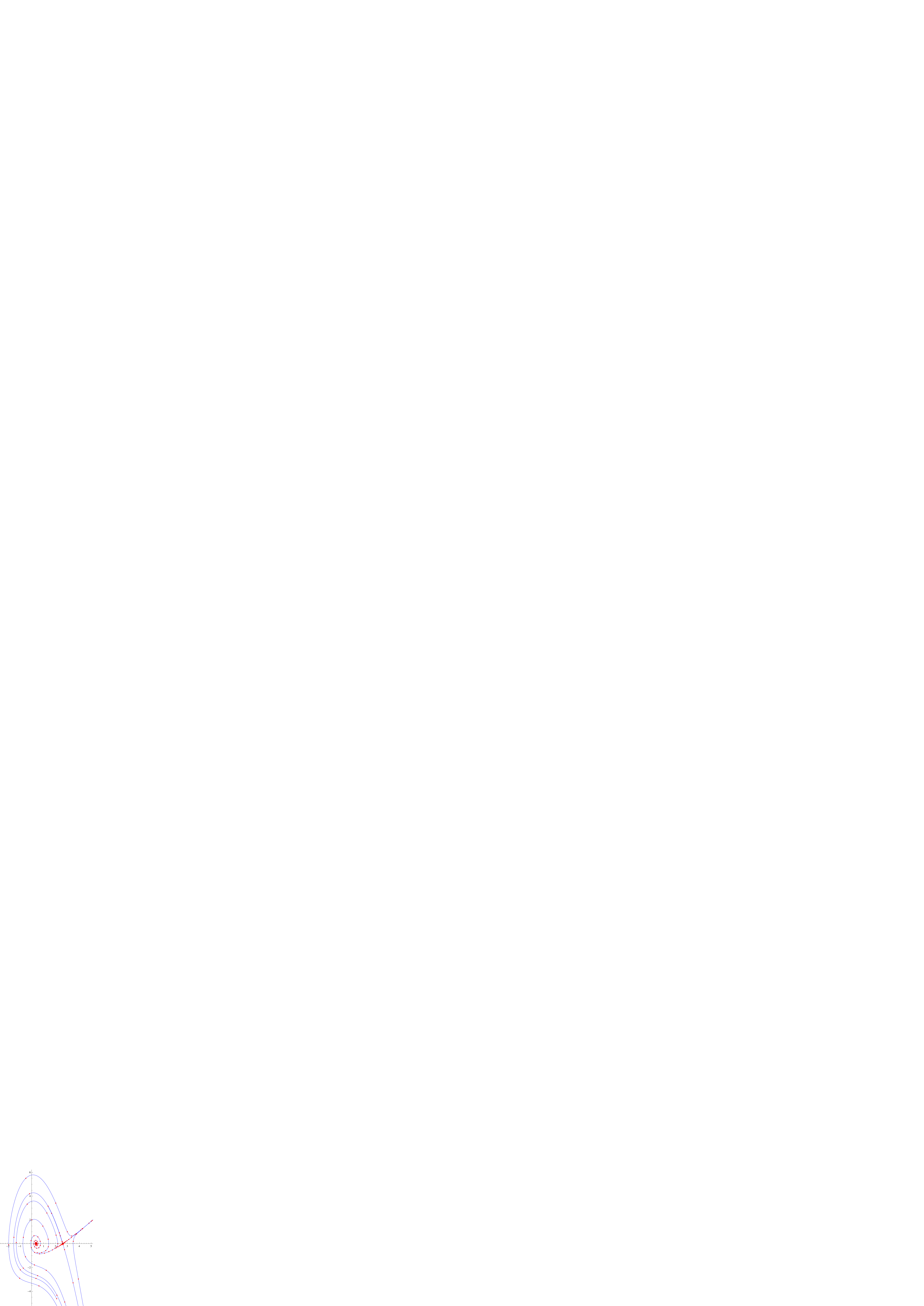}\\ 
\textbf{Figure 3c} curve T-, $\dim_{B}S=1/2$\hskip 1.4cm \textbf{Figure 3d} region $2$ from Figure 2
\end{center}

\begin{center}
\includegraphics[width=4cm,bb=0in 0in 4in 4in]{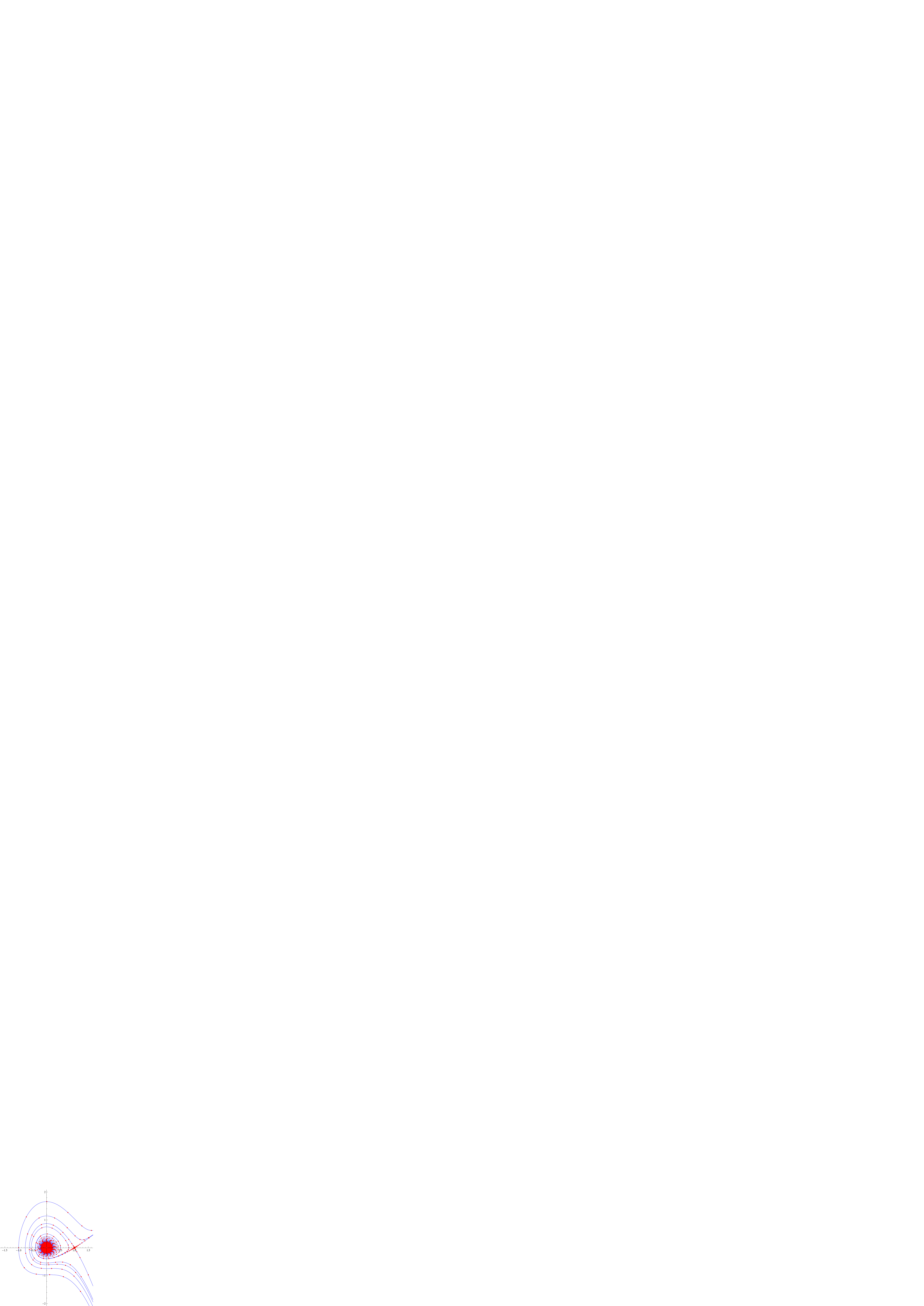} \hskip 2cm \includegraphics[width=2.5cm,bb=0in 0in 4in 4in]{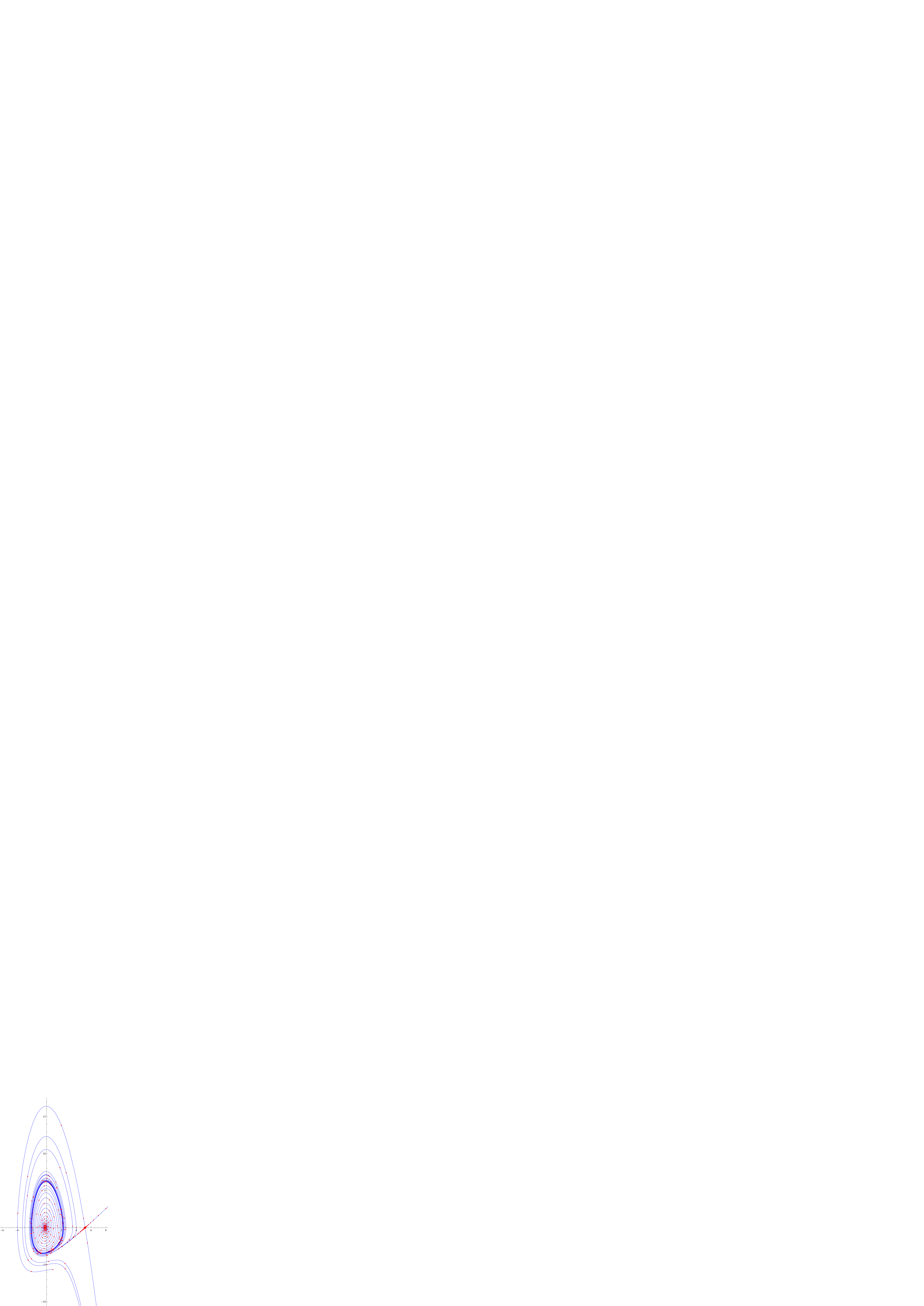}\\ 
\textbf{Figure 3e} curve H, $\dim_{B}S=4/3$ \hskip 1.6cm \textbf{Figure 3f} region $3$ from Figure 2
\end{center}

\begin{center}
\includegraphics[width=3cm,bb=0in 0in 5in 5in]{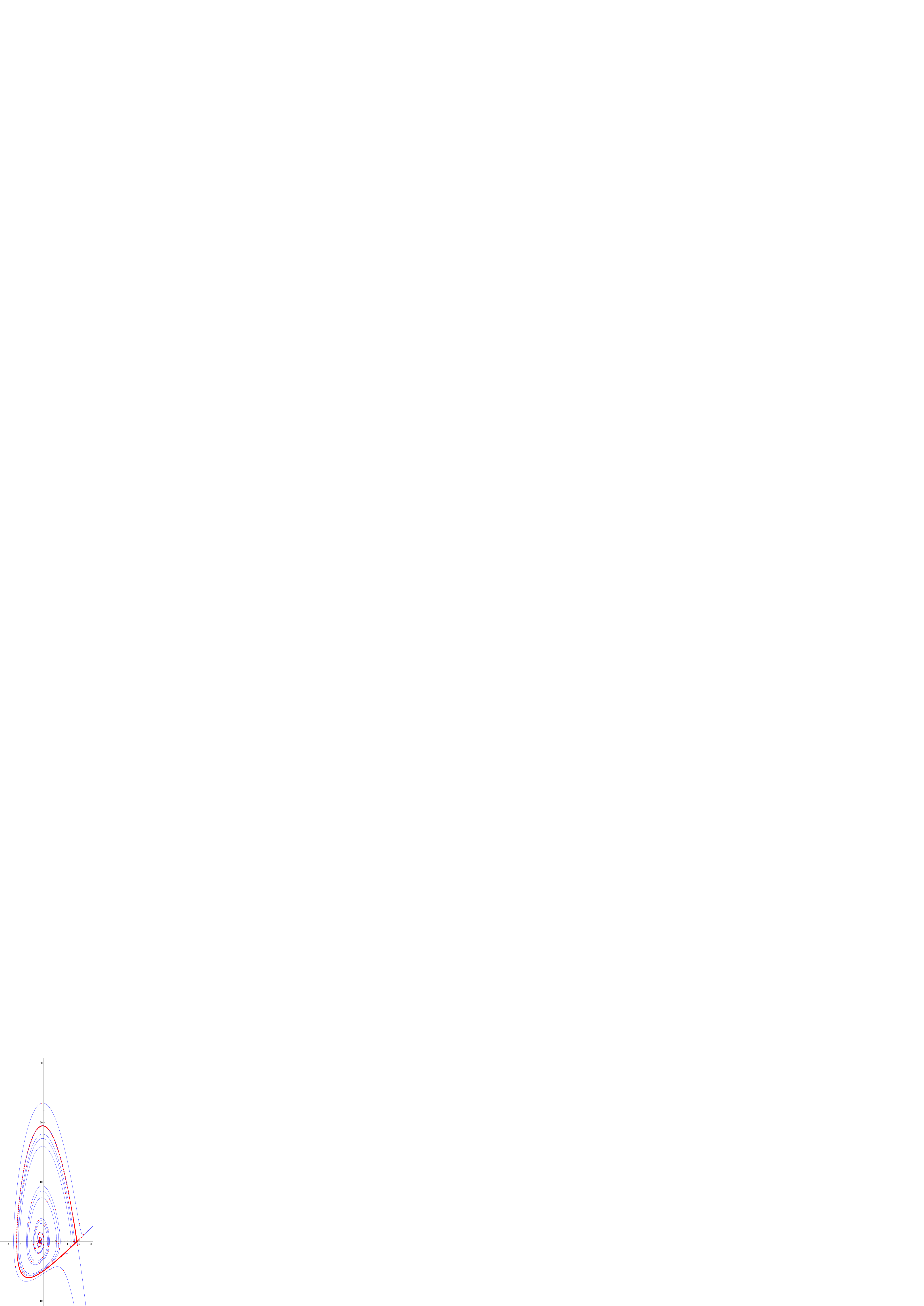} \hskip 2cm \includegraphics[width=6cm,bb=0in 0in 5in 5in]{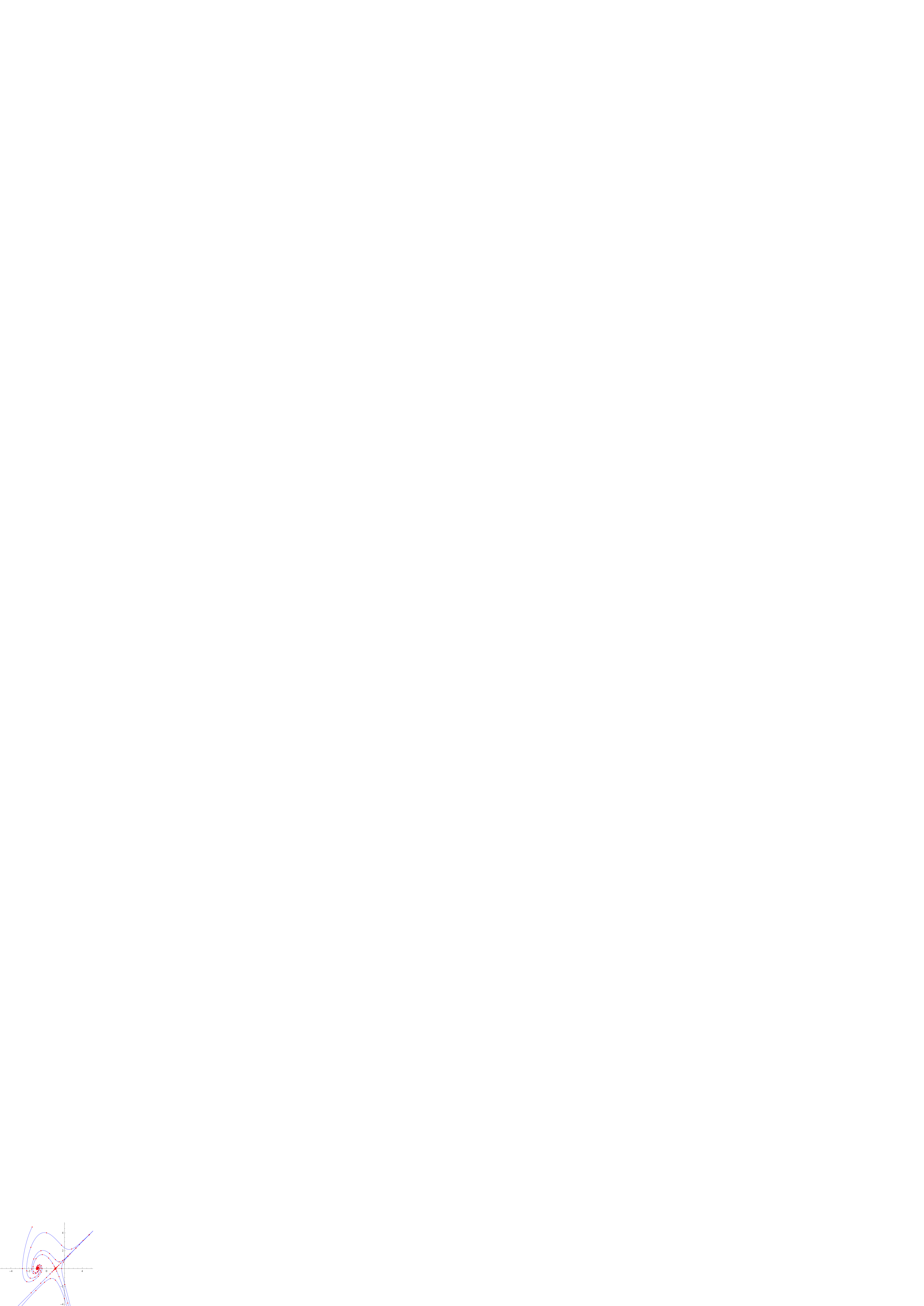}\\ 
\textbf{Figure 3g} curve P \hskip 3cm \textbf{Figure 3h} region $4$ from Figure 2
\end{center}

\begin{center}
\includegraphics[width=3cm,bb=0in 0in 3in 3in]{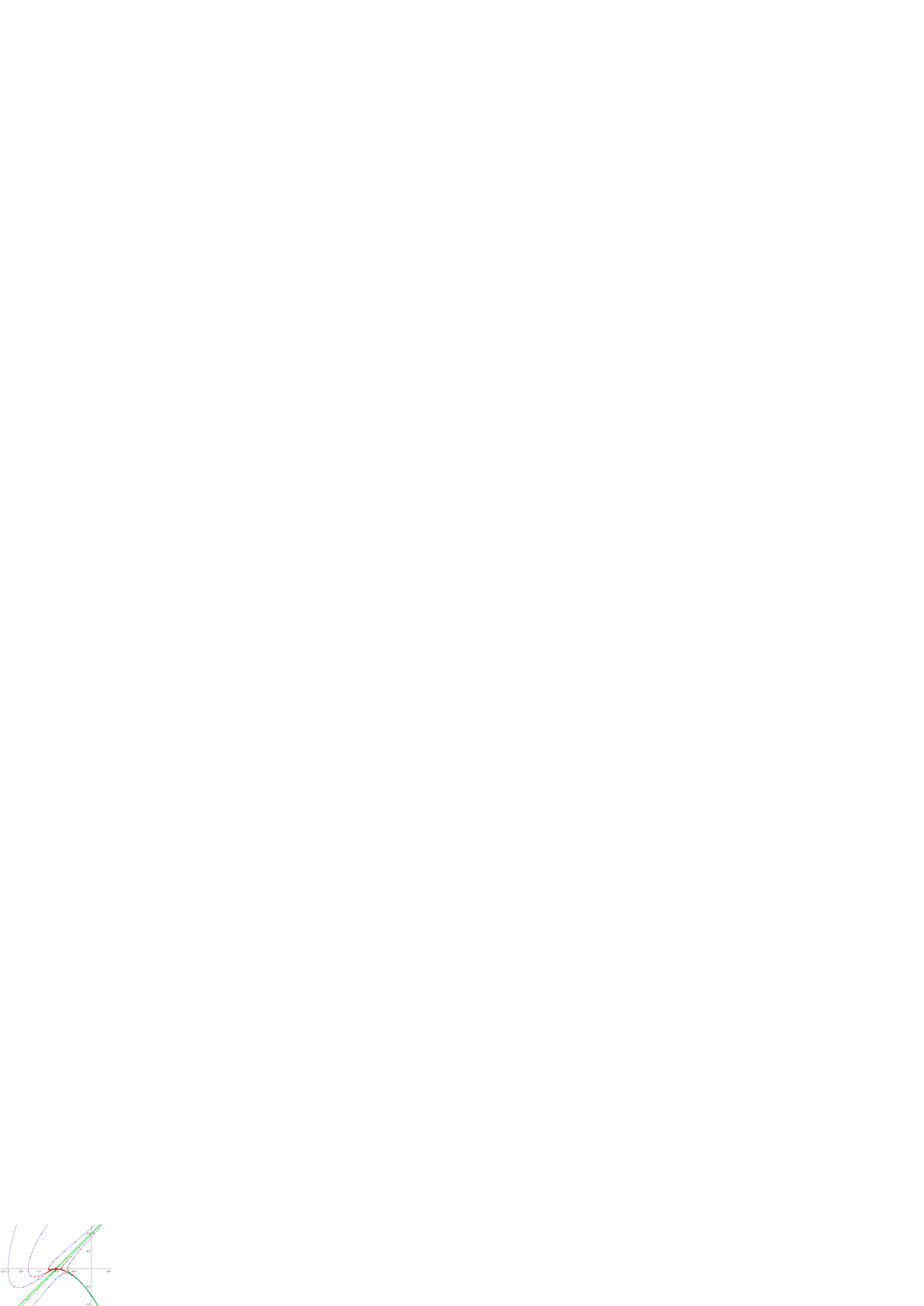} \\ 
\textbf{Figure 3i} curve T+, $\dim_{B}S=1/2$
\end{center}

\subsection{Nilpotent saddle}

The normal form for the two parameter bifurcation of nilpotent saddle is
\begin{eqnarray}
\dot{x}&=&y\nonumber\\
\dot{y}&=& \beta_1 x+\beta_2 y+ x^3-x^2y,
\end{eqnarray}
where $\beta_{1,2}\in\mathbb{R}$ are parameters.  We can see the bifurcation diagram for the degenerate Bogdanov-Takens bifurcation at Figure 4, see more details in \cite{kuz}.
For $\beta_1=\beta_2=0$ we use Theorem \ref{box}, case (1), (i), with $m=3$, $n=2$, $\gamma=\frac{m+1}2=2,$ we get $\dim_{B}S=\dim_{B}S_{x}=1/2$ and $\dim_{B}S_{y}=\frac{1}{3}$. See Figure 5a.  For other cases we use results from \cite{neveda}, \cite{laho}, and \cite{laho2}.
\begin{center}
\includegraphics[width=4cm,bb=0in 0in 2in 2in]{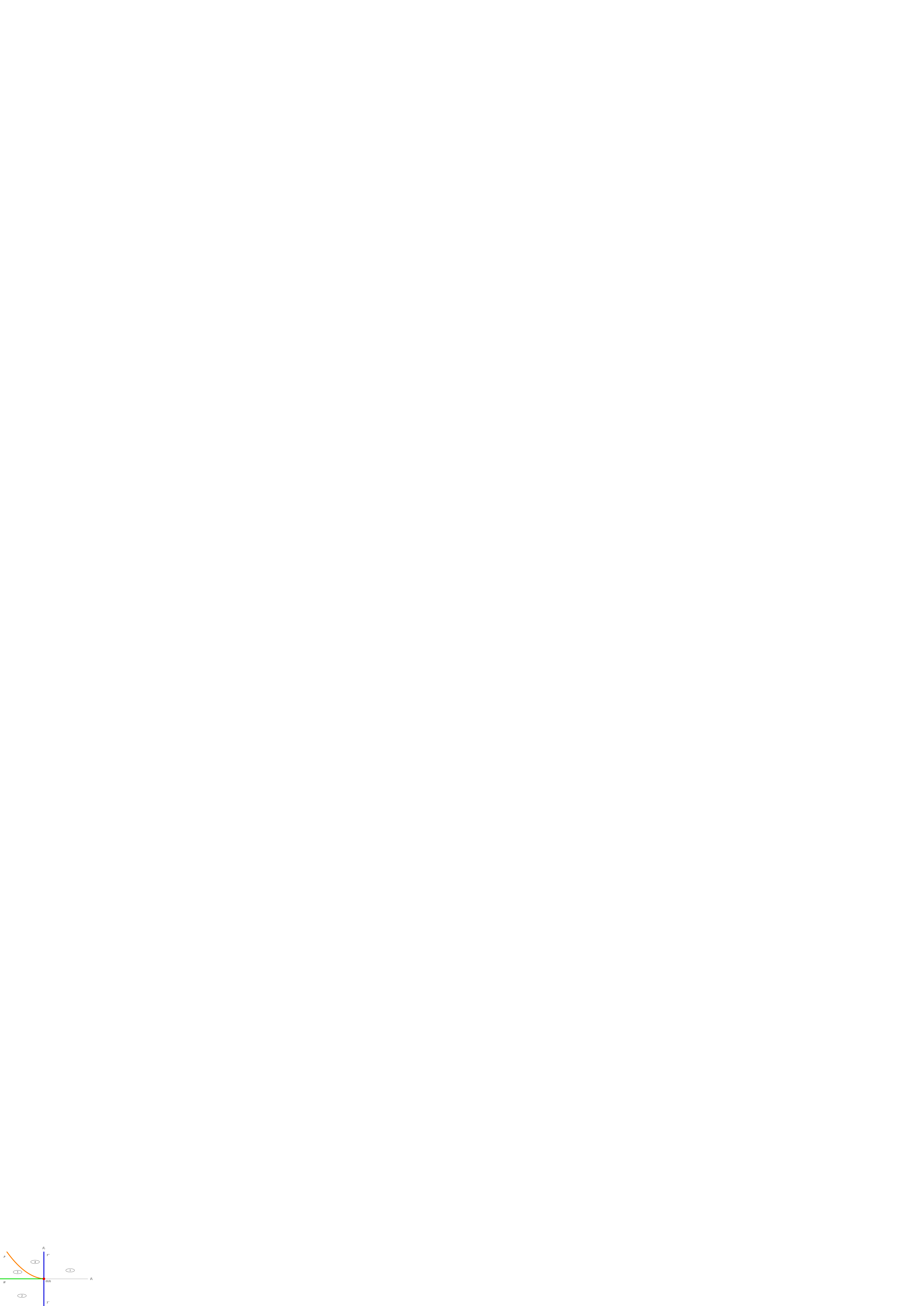} \\ 
 \textbf{Figure 4} bifurcation diagram 
\end{center}

Now through the bifurcation diagram of the unfolding we can see the changing of the box dimension. Similarly as for the cusp case, we start from region 1 with a hyperbolic saddle (Figure 5b), passing through curve $T^{-}$ we get two more singularities. So in region 2 we have two saddles and one node (Figure 5d). This bifurcation is seen by changed box dimension on the curve $T^{-}$ because on one separatrix (orange) is $\dim_{B}S=\dim_{B}S_{x}=\frac{2}{3}$ and $\dim_{B}S_{y}=0$ (Figure 5c). Then on the curve $H$ Hopf bifurcation occurs (Figure 5e), while on $P$ two-saddle loop appears (Figure 5g). The conclusion is analogous as in the previous case, the box dimension is changed at the bifurcation point. Notice that, in this case, the box dimension for $\beta_1=\beta_2=0$ is bigger. It could be connected to the fact that more object is "hidden" in the nilpotent saddle (3 singularities and 1 limit cycle), comparing with the cusp case (2 singularities and 1 limit cycle).

\begin{center}
\includegraphics[width=4cm,bb=0in 0in 4in 4in]{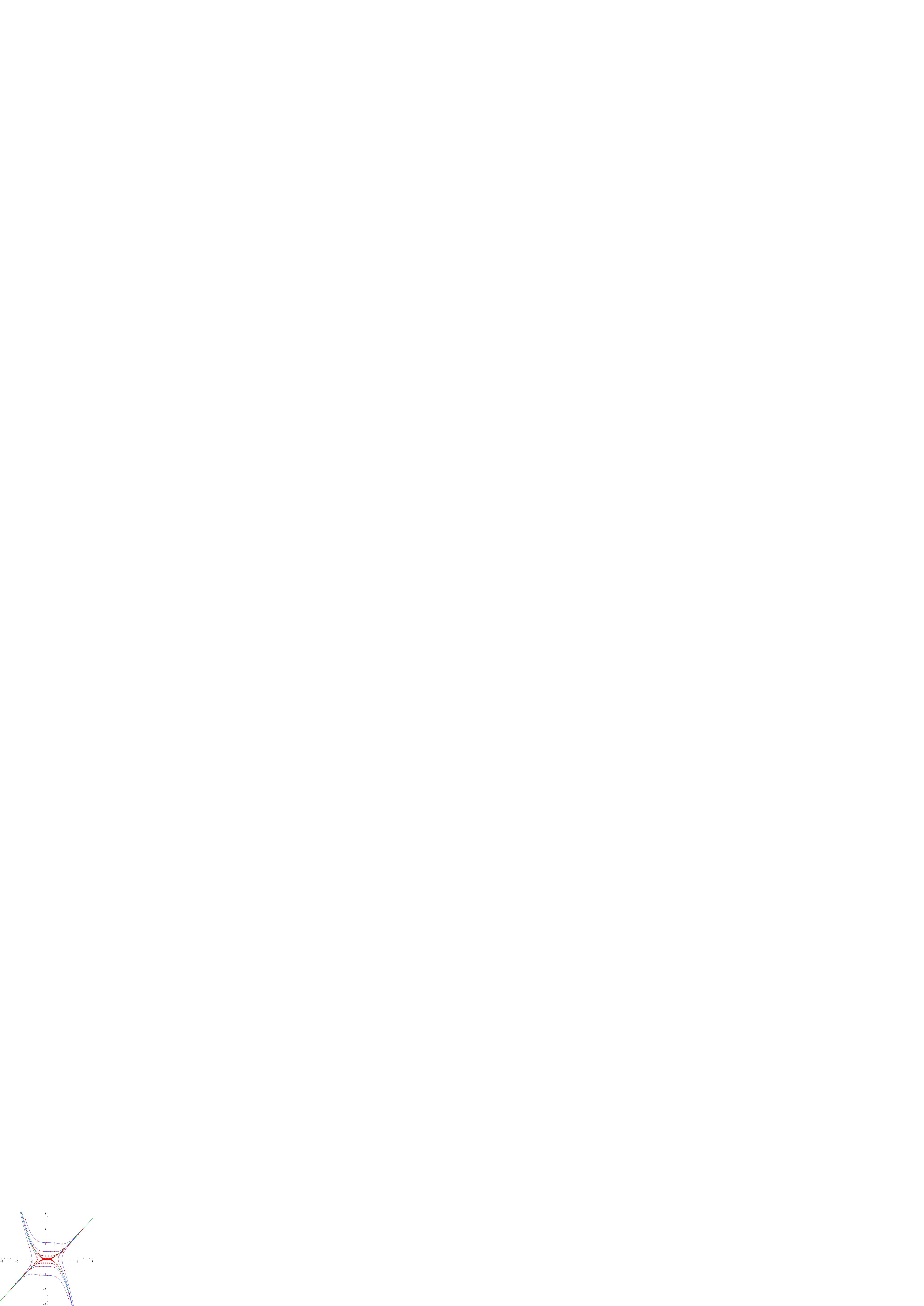} \hskip 2cm \includegraphics[width=4cm,bb=0in 0in 4in 4in]{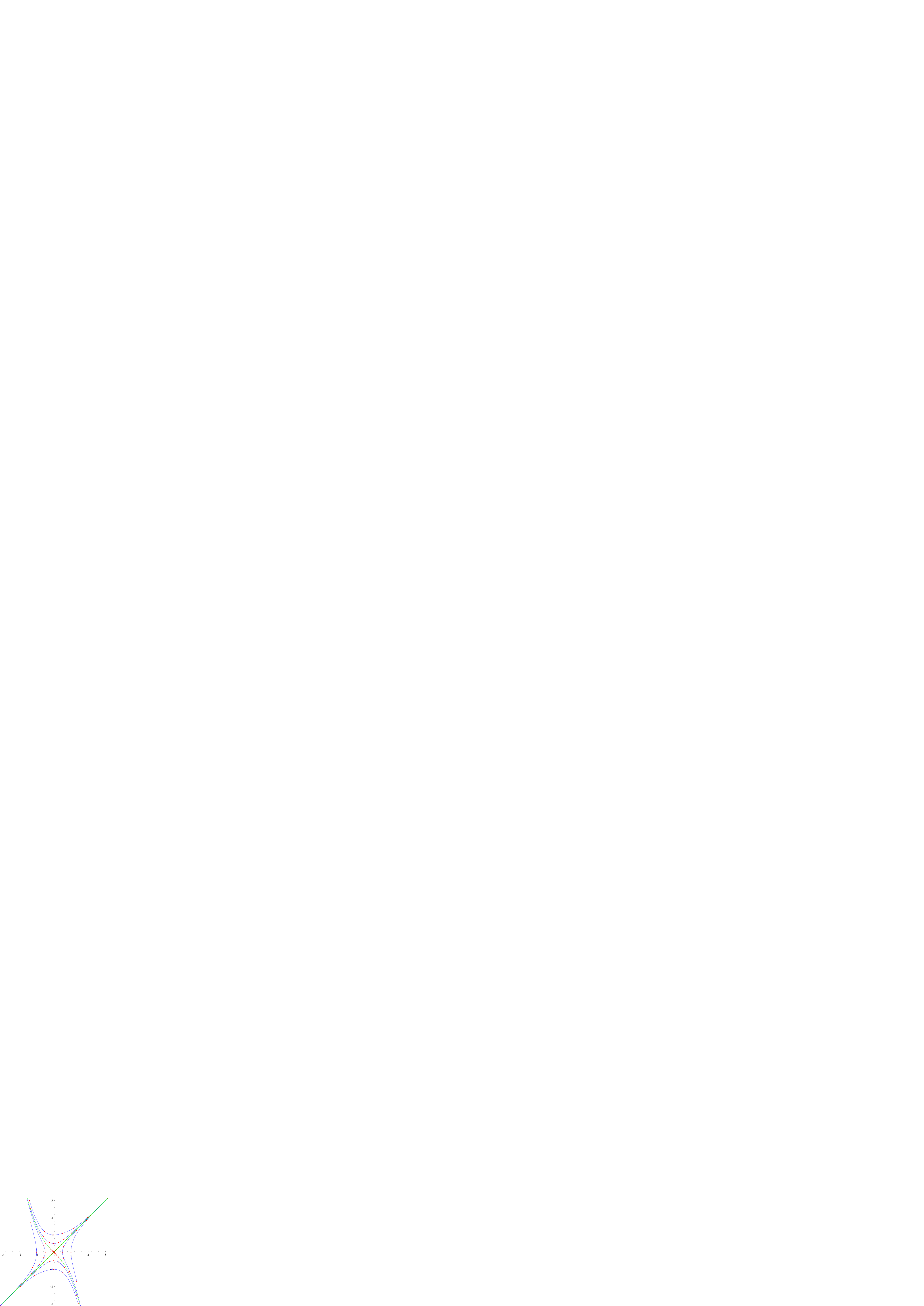}\\ 
\textbf{Figure 5a} $\beta_1=\beta_2=0,\dim_{B}S=1/2$\hskip 1cm \textbf{Figure 5b} region $1$ on Figure 4
\end{center}

\begin{center}
\includegraphics[width=4cm,bb=0in 0in 4in 4in]{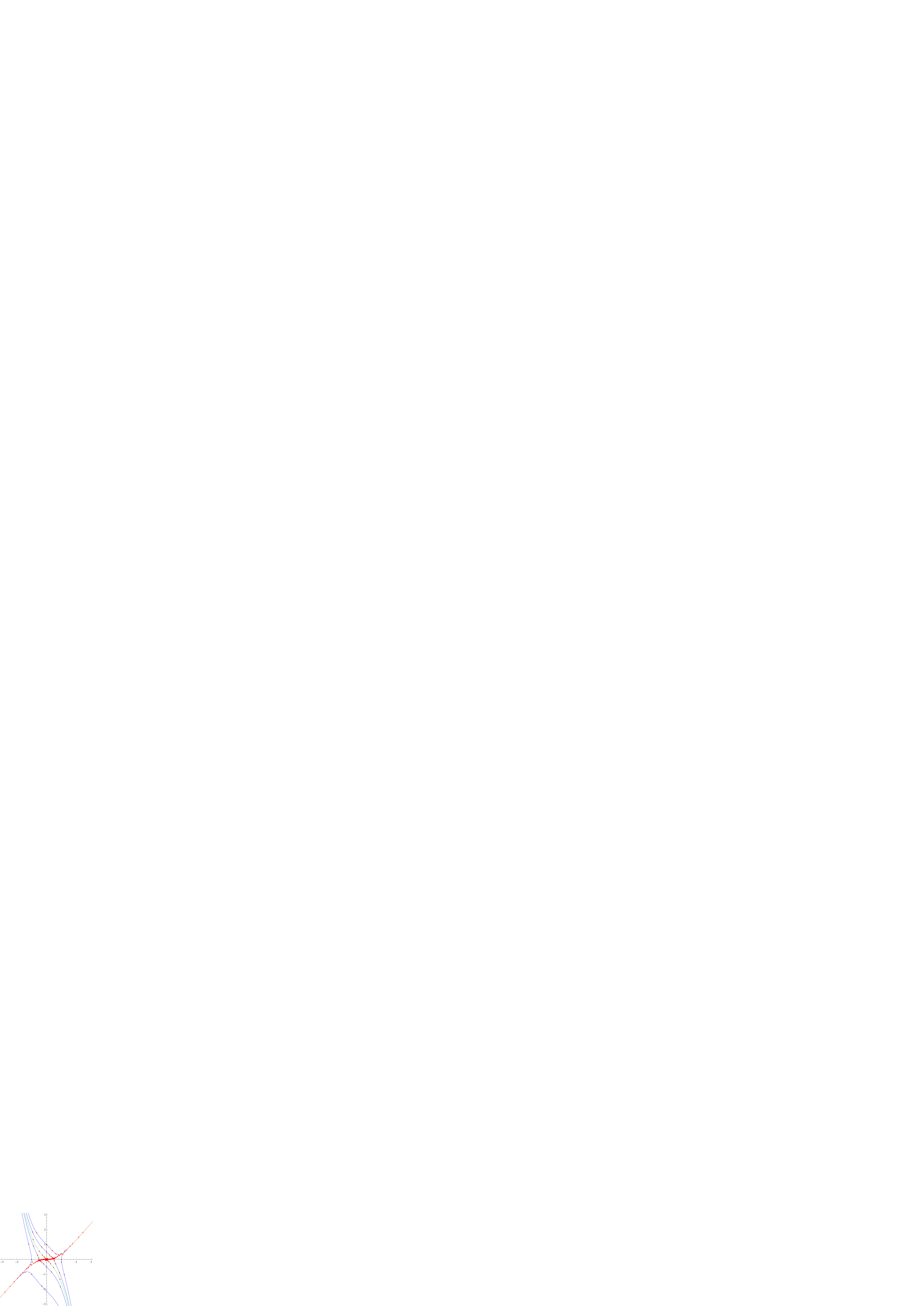} \hskip 2cm \includegraphics[height=4cm,bb=0in 0in 4in 4in]{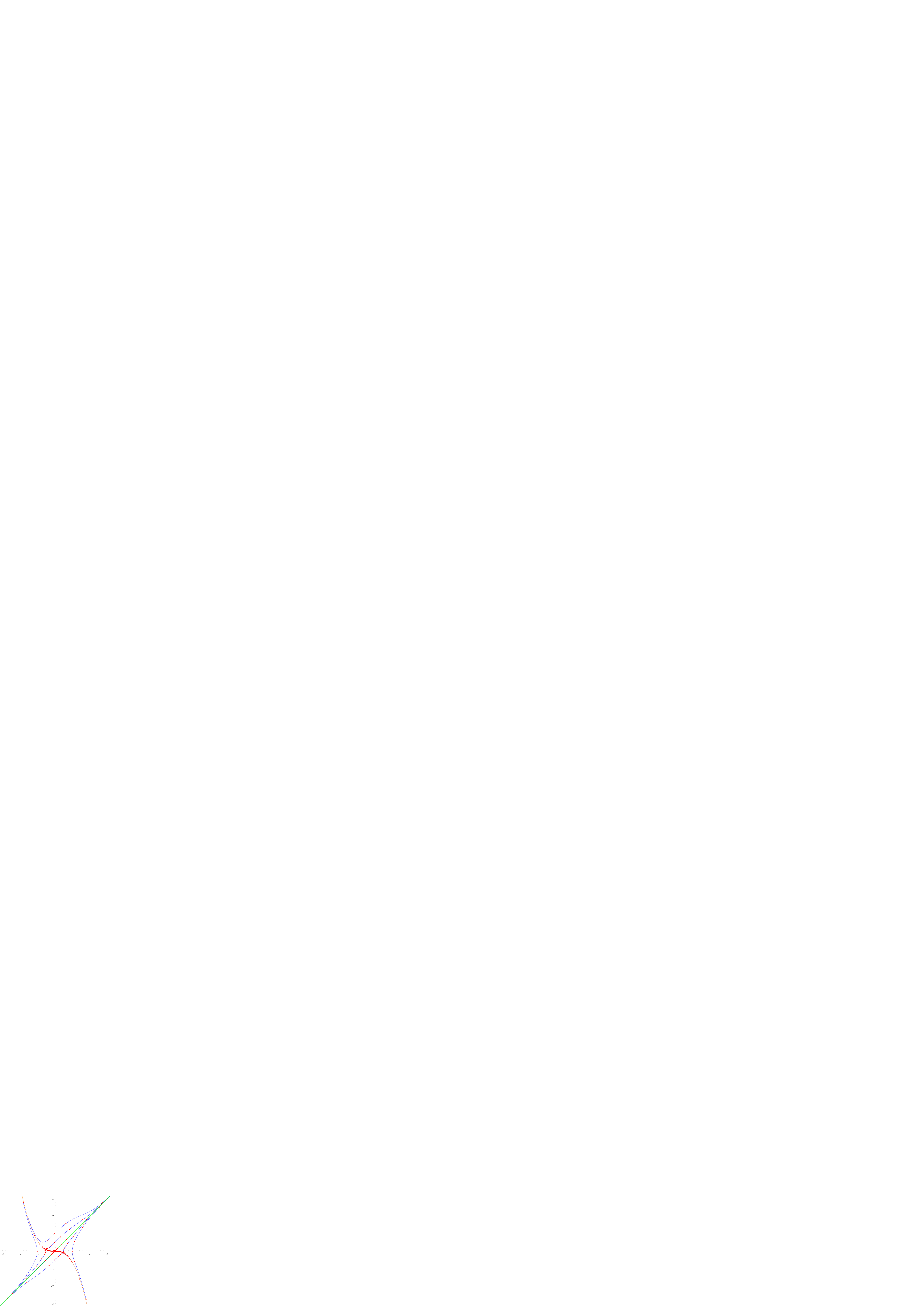}\\ 
\textbf{Figure 5c} curve T-, $\dim_{B}S=2/3$\hskip 1.7cm \textbf{Figure 5d} region $2$ on Figure 4
\end{center}

\begin{center}
\includegraphics[width=4cm,bb=0in 0in 4in 4in]{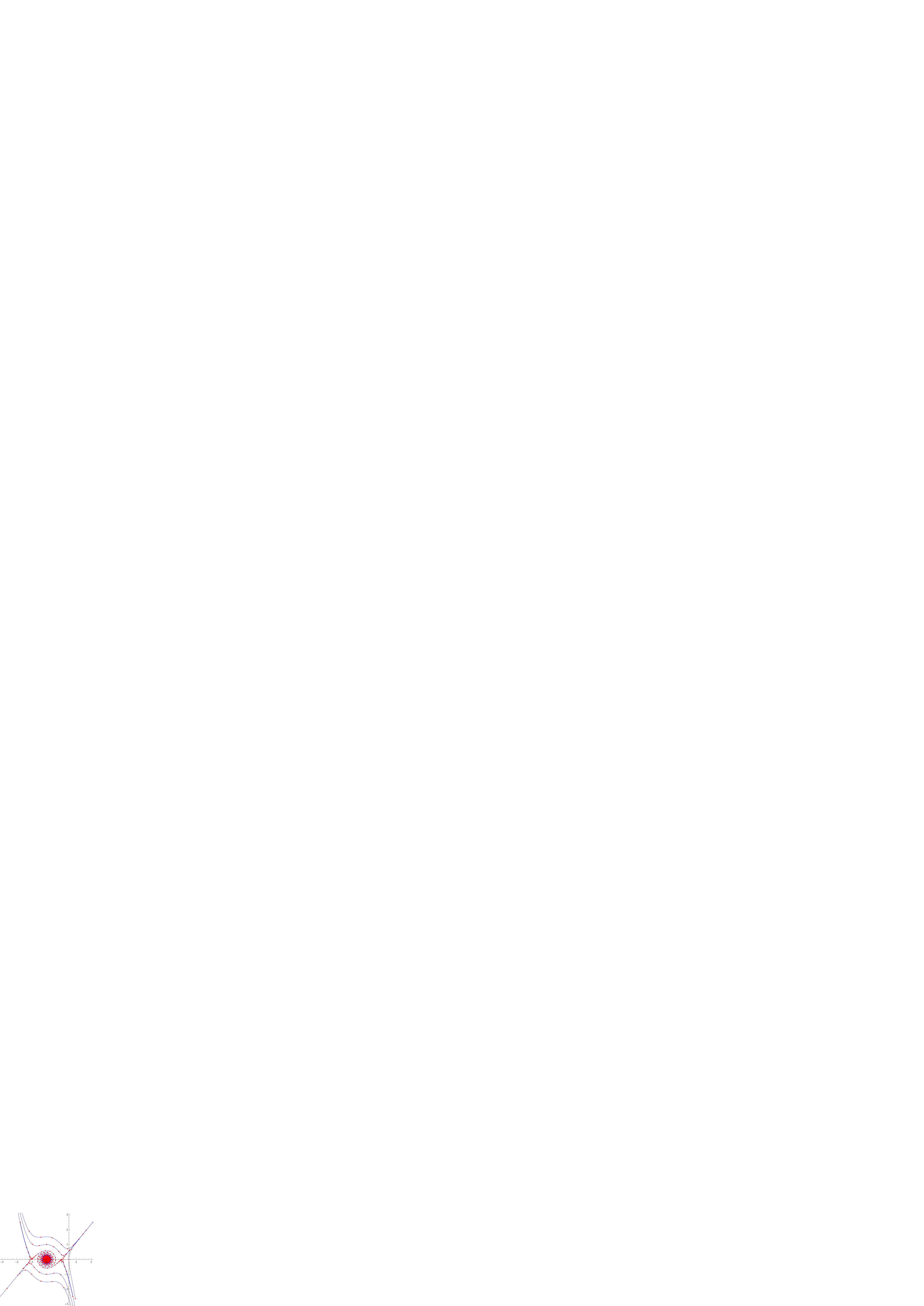} \hskip 2cm \includegraphics[height=4cm,bb=0in 0in 4in 4in]{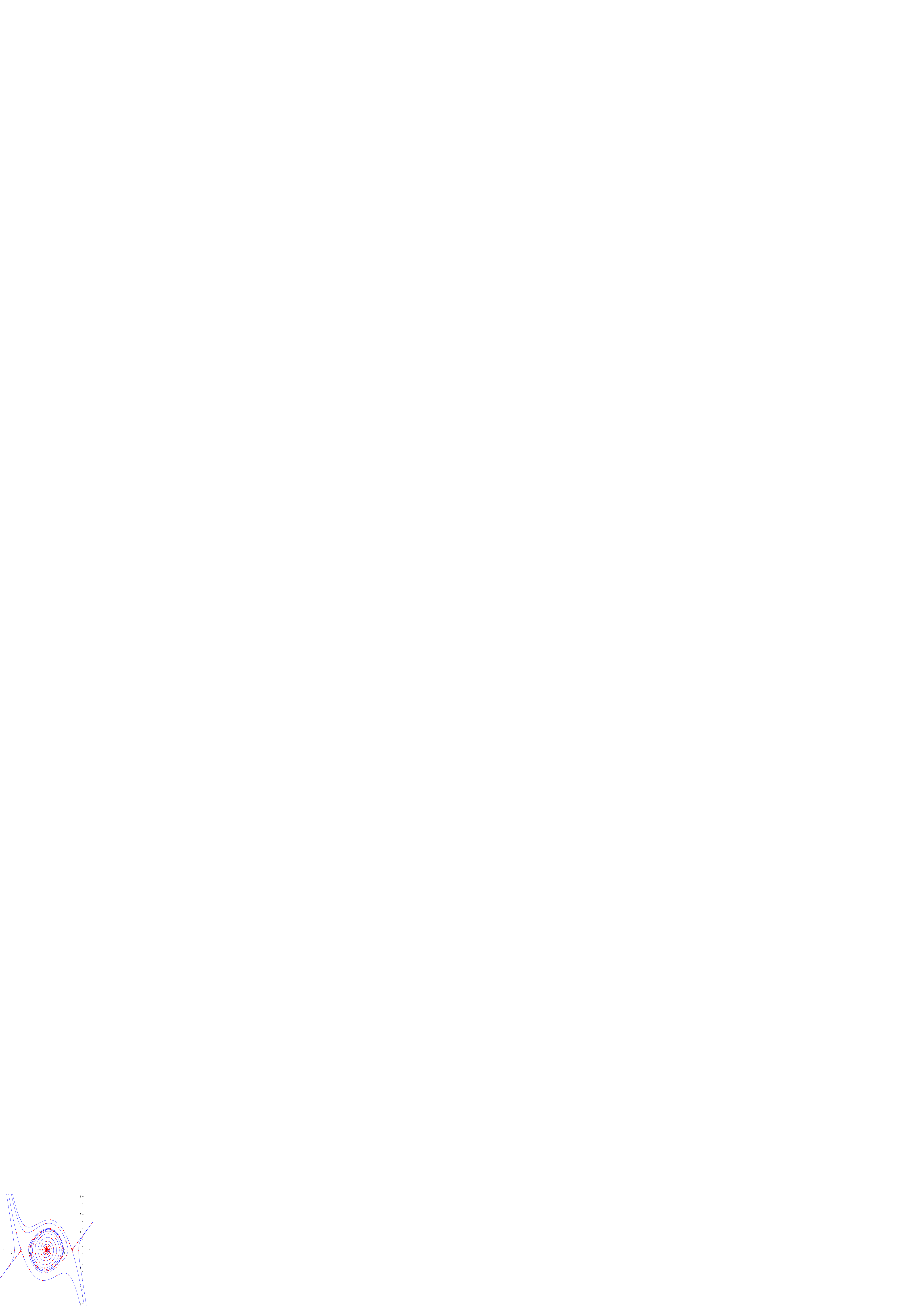}\\ 
\textbf{Figure 5e} curve H, $\dim_{B}S=4/3$\hskip 1.7cm \textbf{Figure 5f} region $3$ on Figure 4
\end{center}

\begin{center}
\includegraphics[height=4cm,bb=0in 0in 4in 4in]{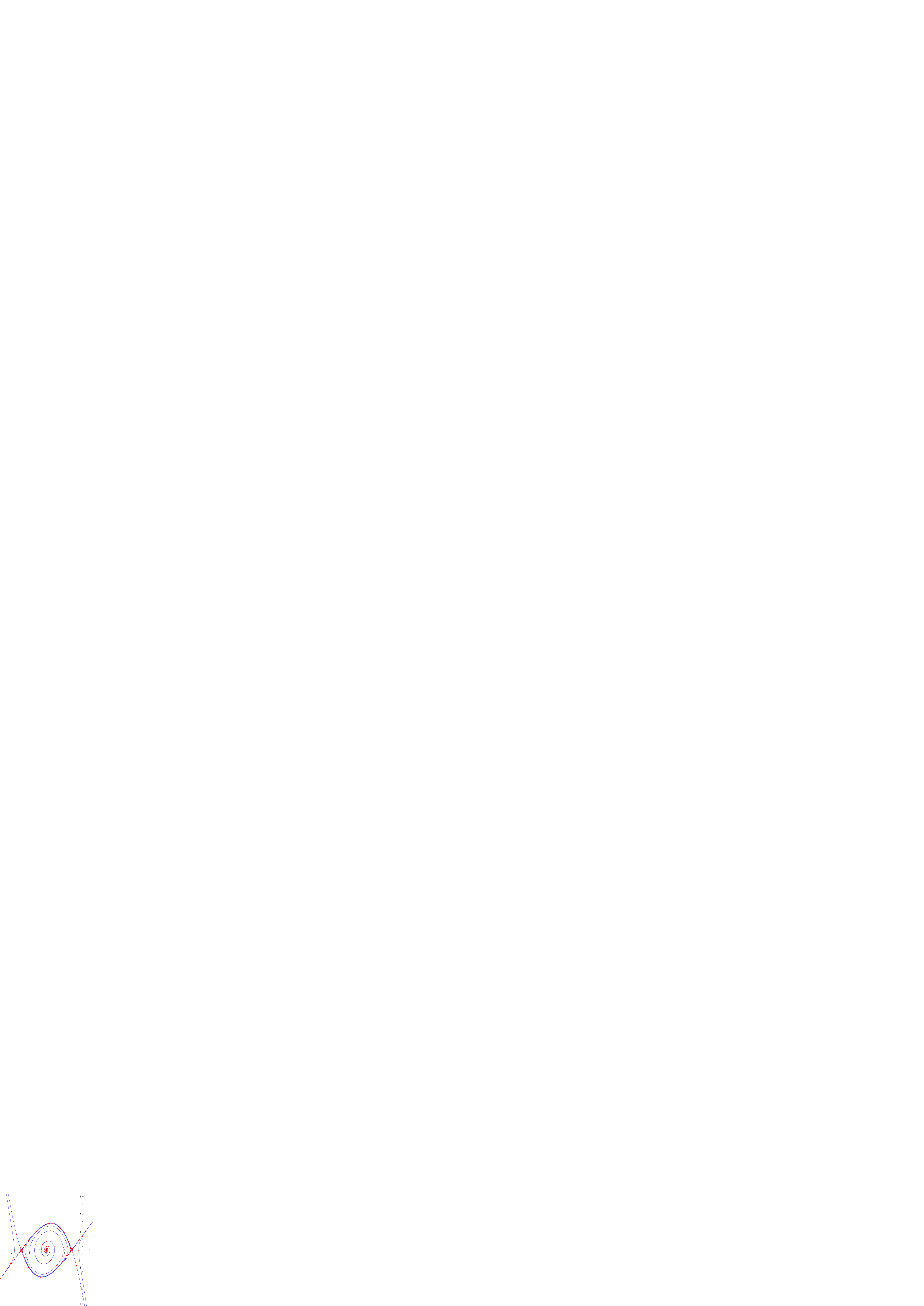} \hskip 2cm \includegraphics[height=4cm,bb=0in 0in 4in 4in]{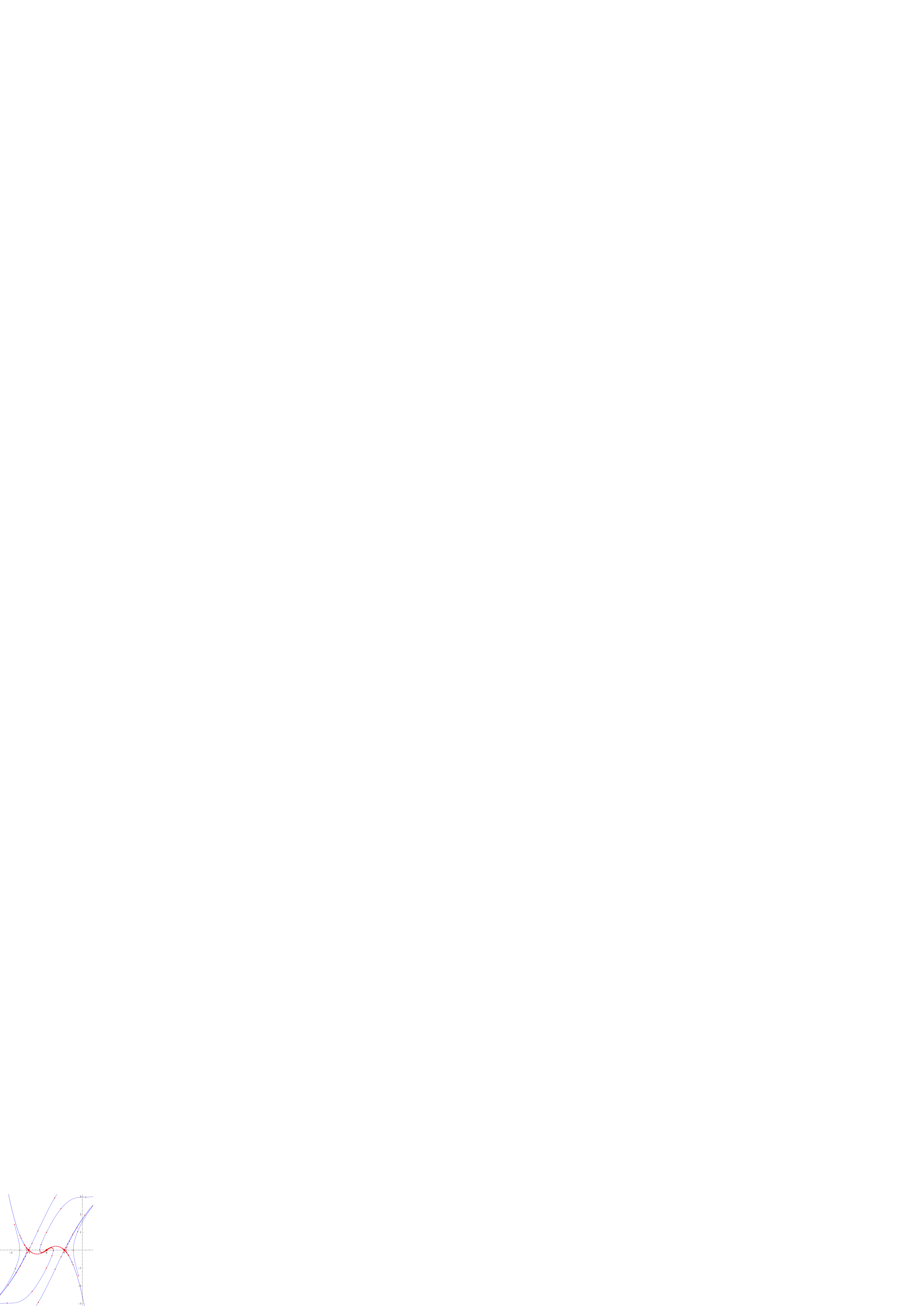}\\ 
\textbf{Figure 5g} curve P \hskip 3cm \textbf{Figure 5h} region $4$ on Figure 4
\end{center}

\begin{center}
\includegraphics[width=4cm,bb=0in 0in 4in 4in]{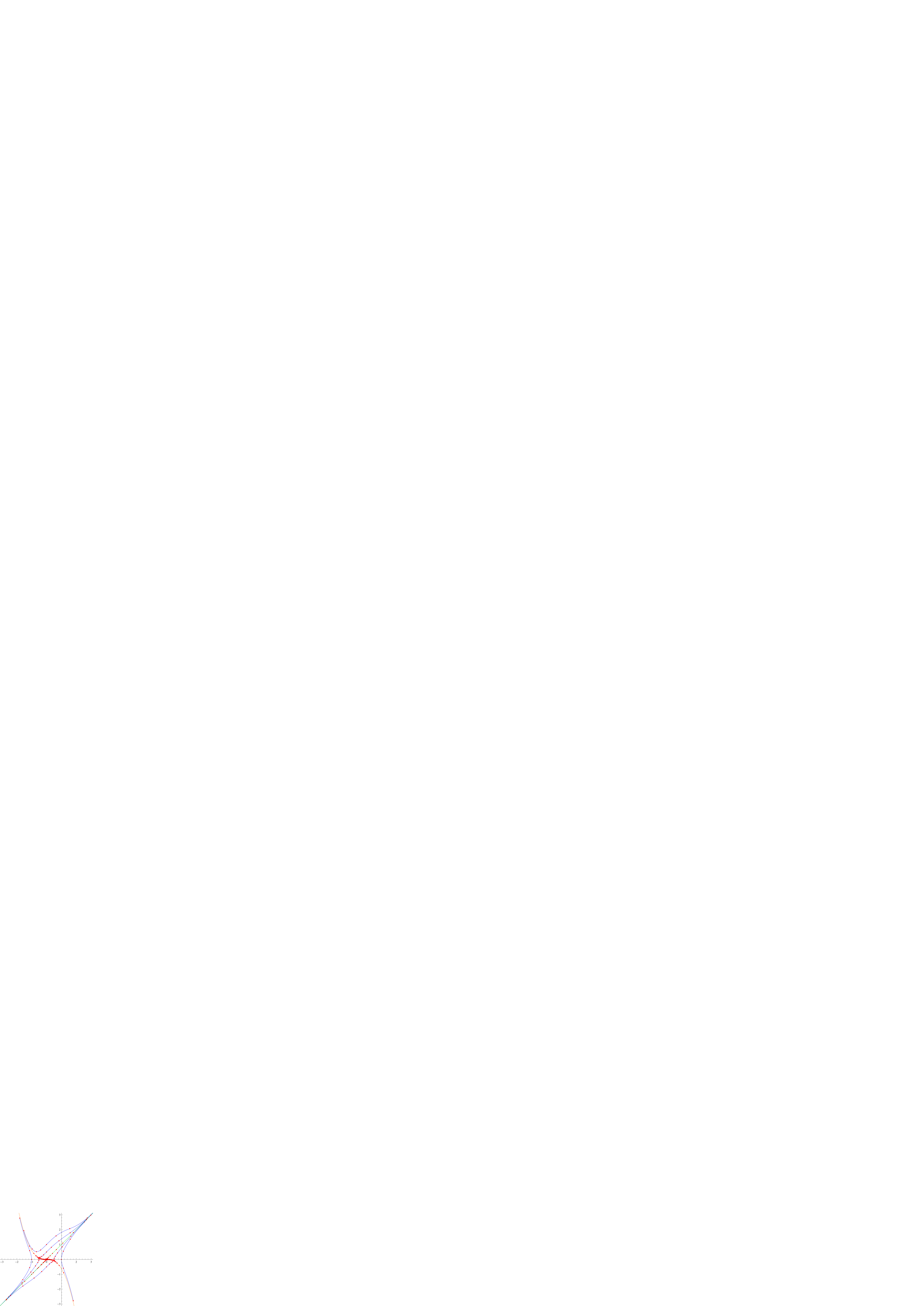} \\ 
\textbf{Figure 5i} curve T+, $\dim_{B}S=2/3$
\end{center}

\subsection{Nilpotent saddle-node}

The system 
\begin{eqnarray}
\dot{x}&=&y\nonumber\\
\dot{y}&=& x^4+xy
\end{eqnarray}
is singular like case $m>2n+1$, $m=4$, $n=1$, from Section 3, and the topological type is a nilpotent saddle-node, see Theorem \ref{type}, (i2).
We use Theorem \ref{box} to obtain the box dimension. For this case, we have behavior $y\simeq x^{\frac{m+1}2}=x^\frac52$ on the separatrix, and using Theorem \ref{box} (2), (i) for $\gamma=5/2$ we obtain $\dim_B S=\dim_{B}S_{x}=1-\frac25=\frac35$ and $\dim_{B}S_{y}=\frac{3}{8}$, where $S$ is a discrete orbit generated by the unit-time map on the separatrix. See Figure 6.

\begin{center}
\includegraphics[width=4cm,bb=0in 0in 3in 3in]{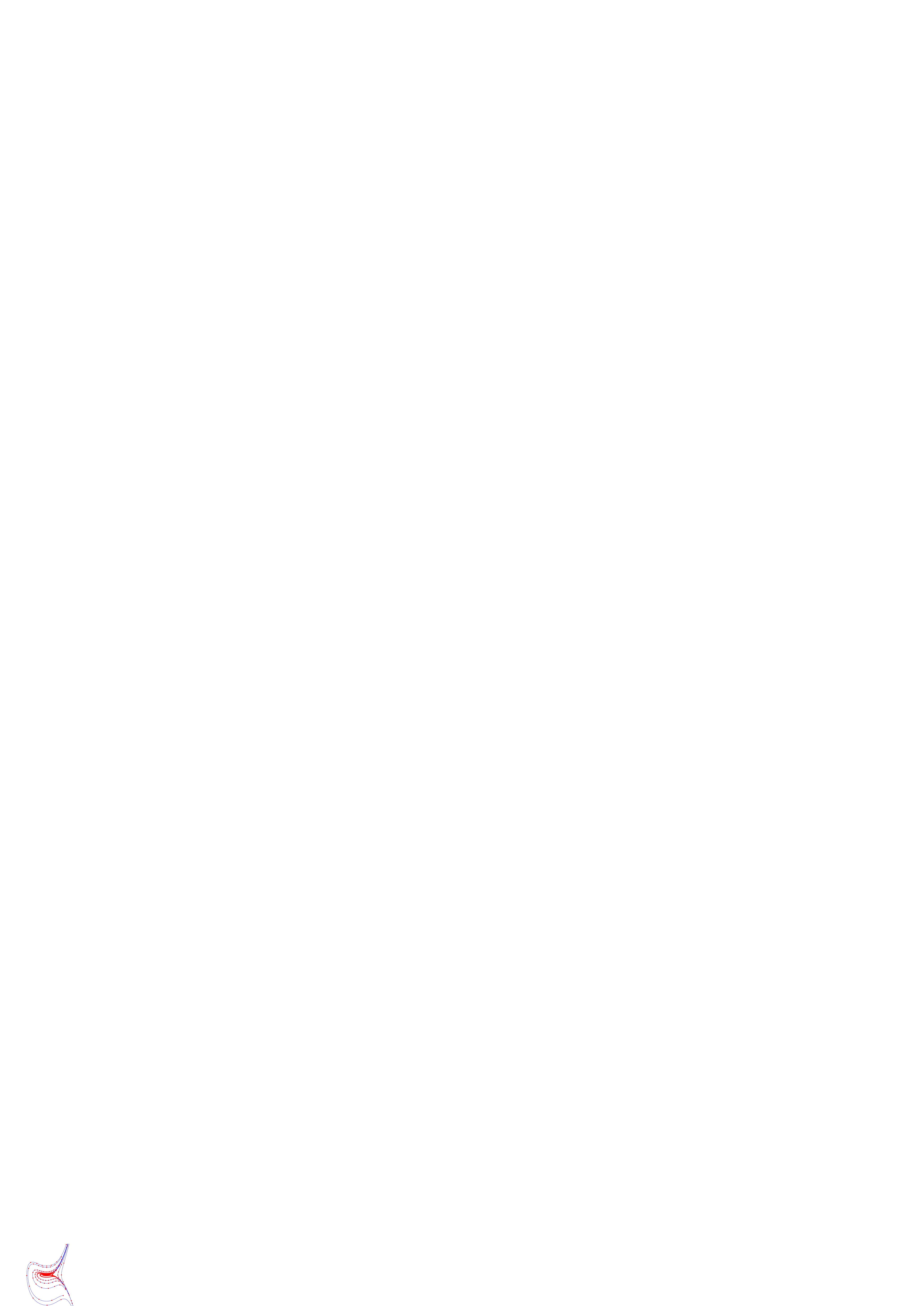} \\ 
\textbf{Figure 6} nilpotent saddle-node, $\dim_{B}S=3/5$
\end{center}

\subsection{Nilpotent node}
The system
\begin{eqnarray}
\dot{x}&=&y\nonumber\\
\dot{y}&=& -x^5-4x^2y
\end{eqnarray}
is mixed case $m=2n+1$, $m=5$, $n=2$, from Section 3, and the topological type is nilpotent node, see Theorem \ref{type}, (iii3).
We use Theorem \ref{box}  to obtain box dimension. For this case,  we have behavior $y\simeq x^{{n+1}}=x^3$   on the separatrix, and using Theorem \ref{box} (1), (i) for $\gamma=3$ we obtain $\dim_B S=\dim_{B}S_{x}=1-\frac13=\frac23$ and $\dim_{B}S_{y}=\frac{2}{5}$, where $S$ is discrete orbit generated by unit-time map on the separatrix. See Figure 7.
\begin{center}
\includegraphics[width=4cm,bb=0in 0in 3in 3in]{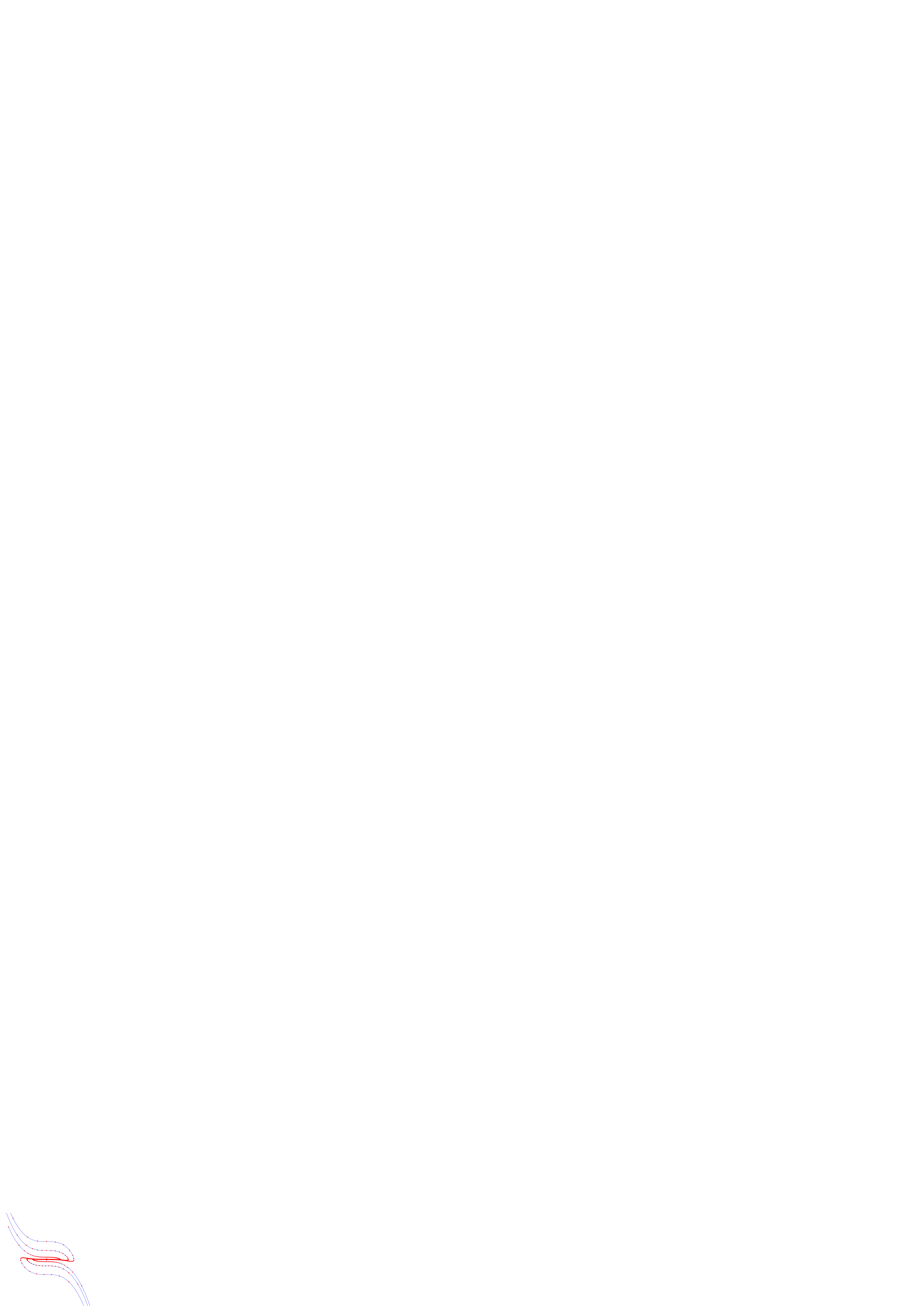} \\ 
\textbf{Figure 7} nilpotent node, $\dim_{B}S=2/3$
\end{center}

\subsection{Hyperbolic and elliptic sector}
The system
\begin{eqnarray}
\dot{x}&=&y\nonumber\\
\dot{y}&=& -x^3+3xy
\end{eqnarray}
is mixed case $m=2n+1$, $m=3$, $n=1$, from Section 3, and the topological type is a singular point with one hyperbolic and one elliptic sector, see Theorem \ref{type}, (iii2).
We use Theorem \ref{box} to obtain the box dimension. For this case, we have the behavior $y\simeq x^{{n+1}}=x^2$ on the separatrix, and using Theorem \ref{box} (1), (i) for $\gamma=2$ we obtain $\dim_B S=\dim_{B}S_{x}=1-\frac12=\frac12$ and $\dim_{B}S_{y}=\frac13$, where $S$ is a discrete orbit generated by the unit-time map on the separatrix. See Figure 8.
\begin{center}
\includegraphics[width=4cm,bb=0in 0in 2.5in 2.5in]
{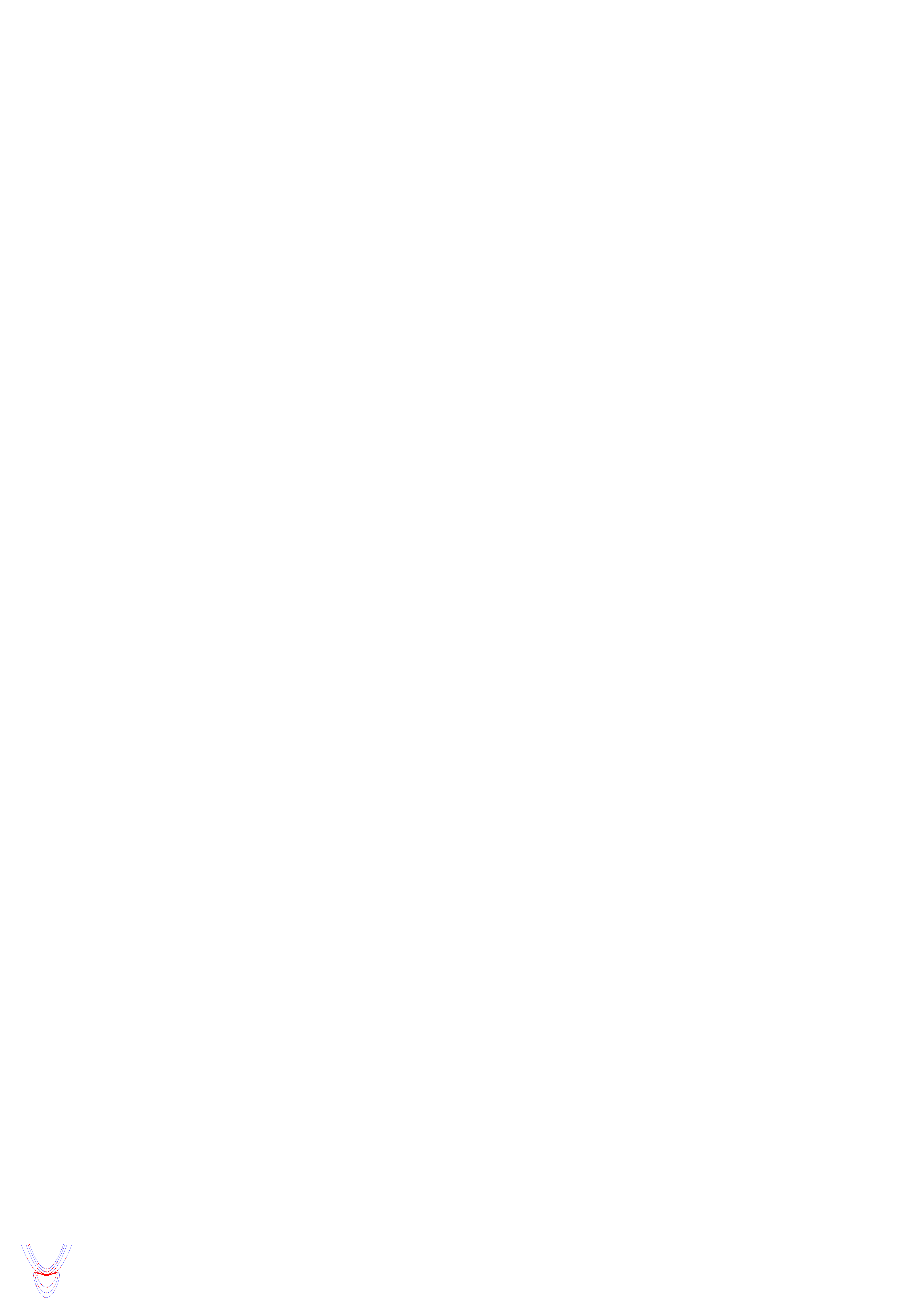} \\ 
\textbf{Figure 8} hyperbolic and elliptic sector, $\dim_{B}S=1/2$
\end{center}

\section{Singularities at infinity for normal form of nilpotent singularity}

The Poincar\'e compactification (see \cite{dla}) is a standard tool for studying the singularities at infinity. In the Poincar\'e compactification, the infinity is presented by a circle.

We continue with the system having a nilpotent singularity at $(0,0)$
\begin{eqnarray} \label{nilsing}
\dot{x}&=&y\nonumber\\
\dot{y}&=&ax^{m}+bx^{n}y, \q a,b\ne 0.
\end{eqnarray}
The case $m<n+1$:\\
Using the formulas from \cite{dla}, p.152, where
$$
x=\frac1v, \q
y=\frac uv,
$$
and dividing by common divisor, in the first chart we obtain
\begin{eqnarray} \label{nilsing1karta}
\dot{u}&=&bu+av^{n+1-m}-u^2v^{n}\nonumber\\
\dot{y}&=& -uv^{n+1}.
\end{eqnarray}
Box dimension of discrete system generated by the unit-time map on the center manifold $v\simeq u^{n+1-m}$ is $\dim_{B}S(x_1,y_1)=1-\frac{1}{2n-m+2}$, where $(x_1,y_1)$ is a first iteration.

Figures 9 represent the orbits generated by unit-time map of the system ($\ref{nilsing1karta}$) with $m=n=2$ (Figure 9a) and $m=2$, $n=3$ (Figure 9b).
\begin{center}
\includegraphics[width=4cm,bb=0in 0in 3in 3in]{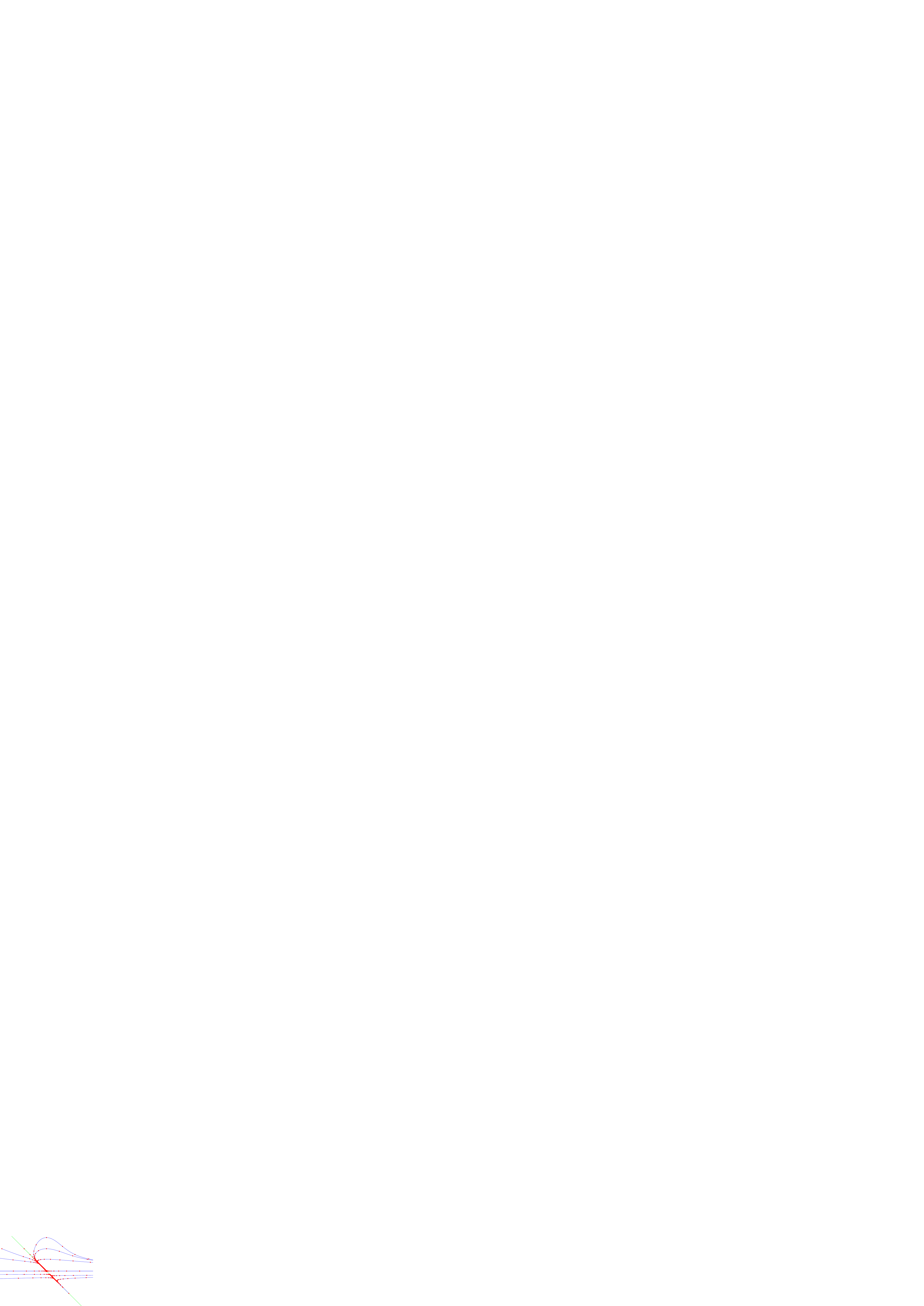}\\
\textbf{Figure 9a} $\dim_{B}S=\frac34$ \\
\includegraphics[width=4cm,bb=0in 0in 3in 3in]{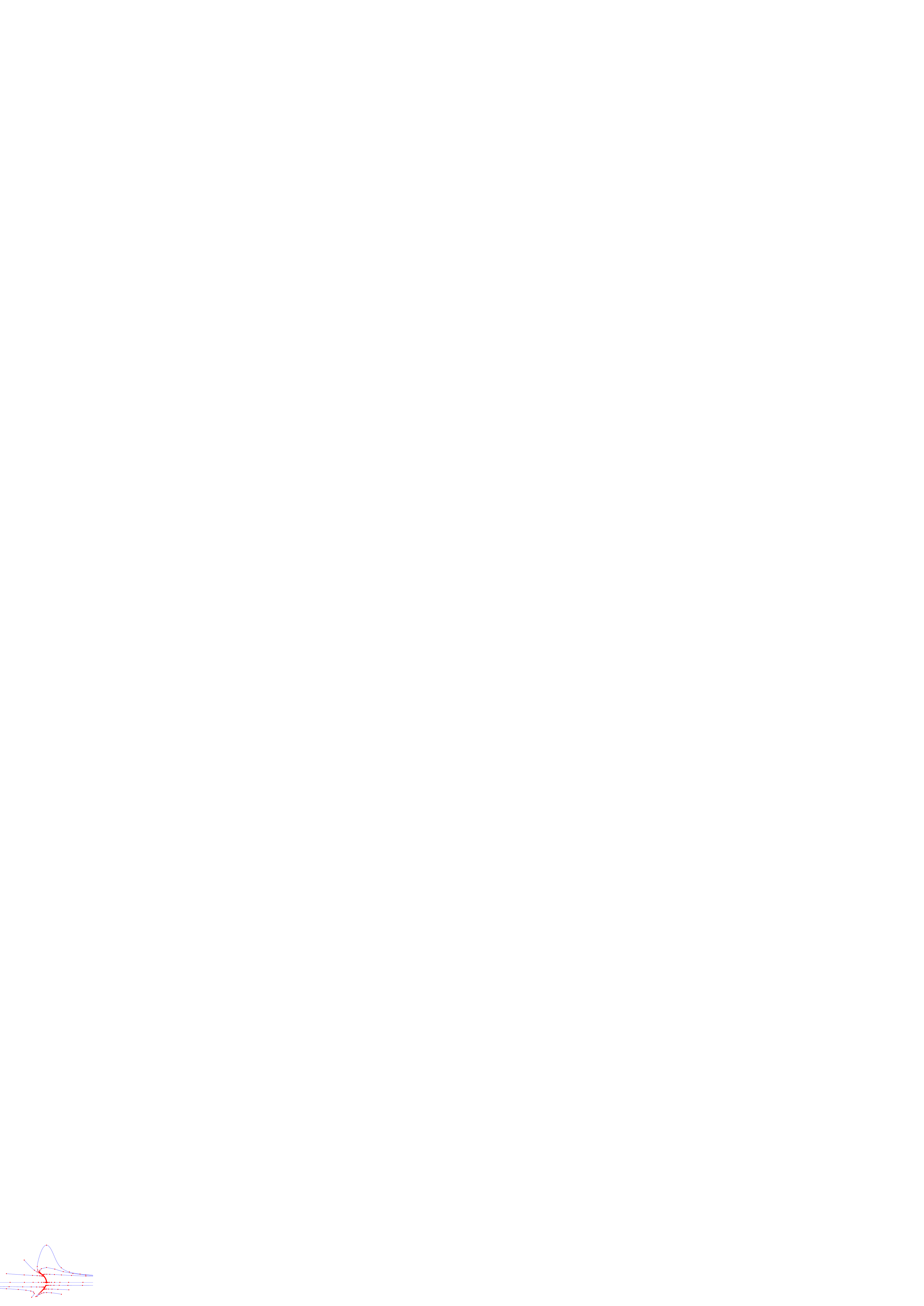}\\
\textbf{Figure 9b} $\dim_{B}S=\frac56$
\end{center}
Using formulas from \cite{dla}, p.152, where
$$
x=\frac uv, \q
y=\frac 1v,
$$
and dividing by common divisor, in the second chart we obtain
\begin{eqnarray} \label{nilsing2karta}
\dot{u}&=&v^{n}-au^{m+1}v^{n+1-m}+bu^{n+1}\nonumber\\
\dot{y}&=& -au^{m}v^{n-m+2}+bu^{n}v.
\end{eqnarray}
Box dimension of the unit-time map on the separatrix $v\simeq u^{\frac{n+1}{n}}$ is $\dim_{B}S(x_1,y_1)=1-\frac{1}{n+1}$.
The result is obtained by computation of the Picard iterations.
Figure 10 represents the unit-time map of the system ($\ref{nilsing2karta}$) with $m=n=2$ (Figure 10a) and $m=2$, $n=3$ (Figure 10b).
\begin{center}
\includegraphics[width=4.5cm,bb=0in 0in 3.5in 4.5in]{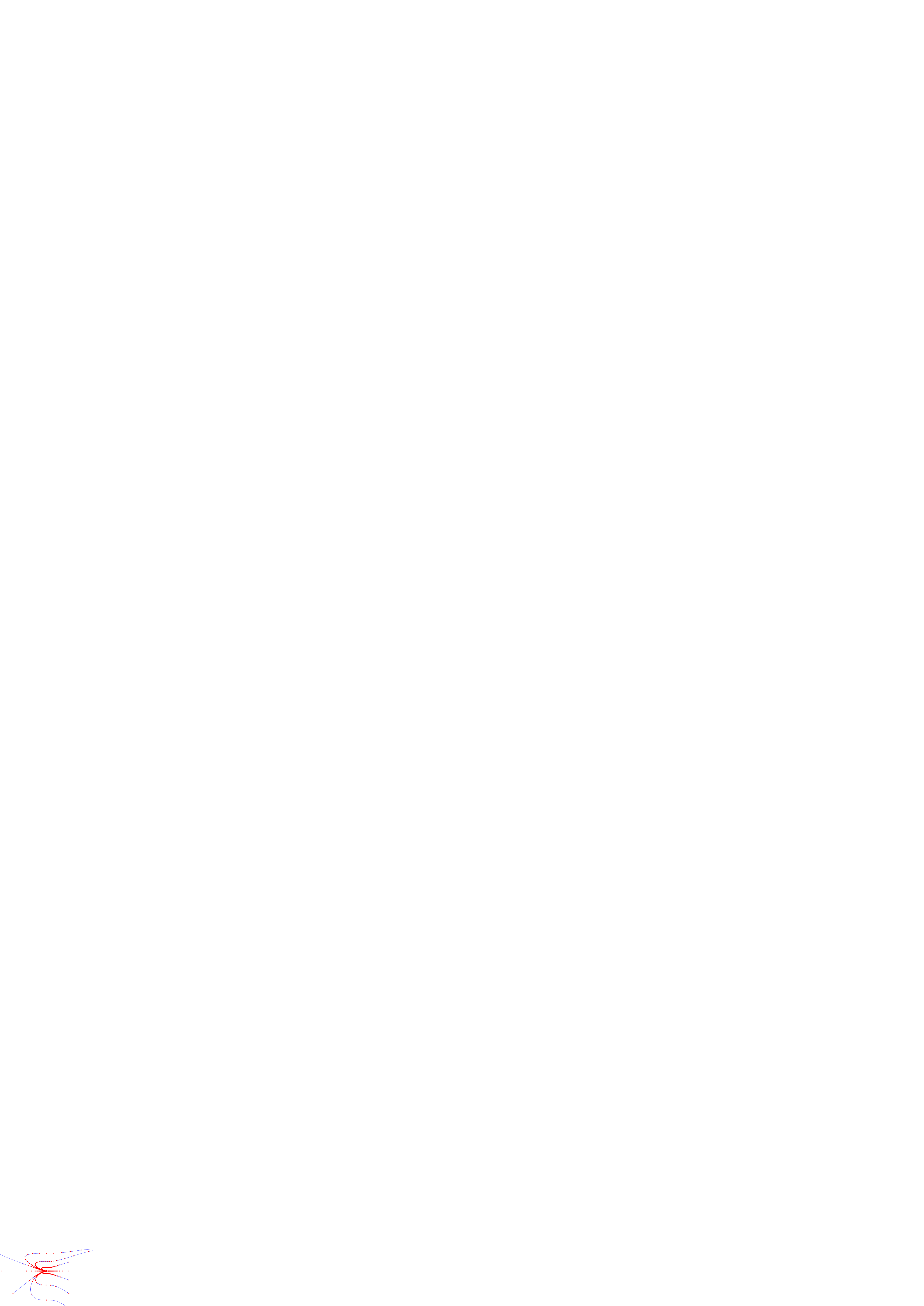}\hskip 2cm\includegraphics[width=4.5cm,bb=0in 0in 3.5in 4.5in]{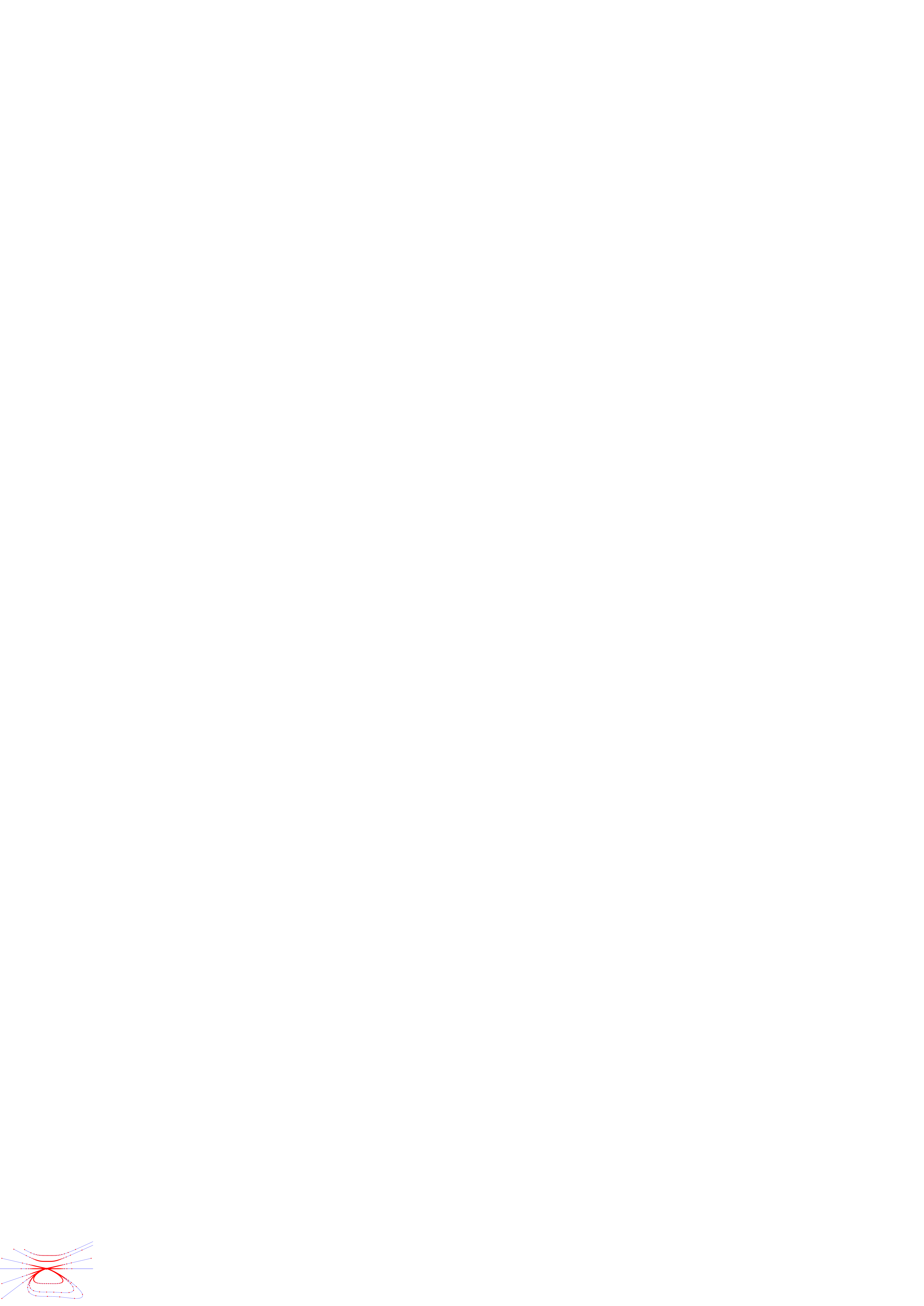}\\  
\textbf{Figure 10a} $\dim_{B}S=\frac23$ \hskip 3cm \textbf{Figure 10b} $\dim_{B}S=\frac34$
\end{center}
The case $m=n+1$ has a singularity in the second chart, and the box dimension is equal as in the previous case.

The case $m>n+1$ has also a singularity in the second chart, and box dimension of the unit-time map on the separatrix $v\simeq u^{\frac{m+1}{m-1}}$ is $\dim_{B}S(x_1,y_1)=1-\frac{1}{m+1}$.

\textbf{Remark 8.} Let us consider the system (\ref{nilsing}) with cusp at origin ($m=n=2$)\\
\begin{eqnarray}
\dot{x}&=&y\nonumber\\
\dot{y}&=&x^{2}+x^{2}y, \q a,b\ne 0.
\end{eqnarray}
Our result show that the box dimension near cusp is $\dim_{B}S=\frac13$, while box dimensions of cusp near infinity in different charts are $\dim_{B}S=\frac34$ and $\dim_{B}S=\frac23$. It can be interesting to find the connection between these values and the cyclicity of cusp.

\section{Singularities at infinity for the normal form of a saddle and dual box dimension}

In this section we study the normal form for a saddle, see e.g. \cite{joyal}, \cite{mardesic}. In  \cite{zuzu} and \cite{belg}, the relation between Lyapunov constants of weak focus and box dimension of spiral trajectory has been established, using the normal form and the Poincar\' e map, respectively. For a saddle there is analogous normal form with coefficients called dual Lyapunov constants or saddle quantities. Here we define dual box dimension using connection between weak focus and saddle, Lyapunov constants and dual Lyapunov constants. 

Since the saddle is hyperbolic, the unit-time map is exponential and we cannot see any interesting behavior of box dimension on the separatrices in the singularity. Situation is different at infinity, we show an example where we can see interesting behavior of box dimension depending on the dual Lyapunov constants. 

We consider the normal form for a saddle
\begin{eqnarray} \label{sedlo}
\dot{x}&=&a_0x+b_0y+\dots +(x^2-y^2)^n(a_nx+b_ny)\\
\dot{y}&=& b_0x+a_0y+\dots +(x^2-y^2)^n(b_nx+a_ny), a_n\ne 0 \nonumber
\end{eqnarray}
where the coefficients $a_k$ are called {\it dual Lyapunov constants or saddle quantities}. If $a_k\ne 0$ is fist nonzero coefficient, we say that a saddle is weak of order $k$, completely analogous to the standard definition of a weak focus of order $k$.  Saddle quantities play important role in the  problem of cyclicity near a saddle loop, they appear in asymptotic expansion of the Dulac map. Poincar\' e map near a saddle loop is a composition of Dulac map (passing near a saddle), and the regular part.

Standard normal form for a weak focus is (see \cite{takens})
\begin{eqnarray} \label{weakfocus}
\dot{x}&=&\hskip-0.2cm-y+a_0x+a_1x(x^2+y^2)\dots +a_n(x^2+y^2)^n\\
\dot{y}&=& x+a_0y+a_1y(x^2+y^2)\dots +a_n(x^2+y^2)^n, a_n\ne 0. \nonumber
\end{eqnarray}
Since $a_n\ne 0$, we write $a_n=1$ and in polar coordinates we have
\bgeqn\label{takenspol}
\dot r&=&r(r^{2n}+\sum_{i=0}^{n-1}a_ir^{2i})\\
\dot\f&=&1.\nonumber
\endeqn
In  \cite{zuzu} we exploited Tricot's formula for the box dimension of a spiral $\Gamma$ defined by $r={\varphi}^{-\alpha}$, $0<\alpha \le1$, $\dim_B\Gamma=\frac{2}{1+\alpha}$. We obtained the result that if $a_k\ne 0$ is a first nonzero coefficient then 
spiral trajectory $\Gamma$ has box dimension $\dim_B\Gamma=\frac{4k}{2k+1}$.

Put $a_n=b_0=1$, $b_i=0$ for $i=1\dots n$ in (\ref{sedlo}), and change of variables 
$$
x=r\cosh\varphi,\q y=r\sinh\varphi
$$
we get exactly the same system as (\ref{takenspol}), but in hyperbolic coordinates.
Solutions of system (\ref{takenspol}) are comparable to $r={\varphi}^{-\alpha}$, $0<\alpha \le1$, with appropriate $\alpha$. If the trajectory $\Gamma_h$ of system (\ref{takenspol}) in hyperbolic coordinates is comparable to
$
r={\varphi}^{-\frac1{2k}},
$ 
we say that $\Gamma_h$ has {\it dual box dimension} $\dim^*_B\Gamma_h=\frac{4k}{2k+1}$.

Let study some example of the system (\ref{sedlo}) at  $\infty$. Assuming that $a_k=b_0=1$, and the other coefficients vanish
we get
\begin{eqnarray} \label{sedloex}
\dot{x}&=&y+(x^2-y^2)^k\\
\dot{y}&=&x+(x^2-y^2)^k. \nonumber
\end{eqnarray}
Using formulas from \cite{dla}
$$
x=\frac1v, \q
y=\frac uv,
$$
and dividing by common divisor, in the first chart we obtain the system with two singularities $(\pm 1, 0)$. We translate it to $(1,0)$, but keep the notation $u,v$, and  get
\begin{eqnarray} \label{nilsedlobes}
\dot{u}&=&-2uv^{2k-1}+u^2v^{2k-1}\nonumber\\
\dot{v}&=& -v^{2k}-uv^{2k}-(-1)^{k}u^{k}(2+u)^{k}.
\end{eqnarray}
The system (\ref{nilsedlobes}) has a singularity $(0,0)$ and the invariant set is $u=0$. Using Picard iterations we get that the unit-time map is comparable to $v-v^{2k}$.  Then the box dimension of an orbit of the unit-time map on the $v$-axes is $\dim_{B}S(v_1)=1-\frac{1}{2k}$, where $(0,v_1)$ is a first iteration. Box dimension increases by $k$, which is the order of the weak saddle.
Results are the same for singularity $(-1,0)$, also the second chart does not show any new singularities.

Figure 11a represents the discrete orbits of the unit-time map of the system ($\ref{nilsedlobes}$) with $k=1$
\begin{eqnarray} \label{nilsedlobes1}
\dot{u}&=&-2uv+u^2v\nonumber\\
\dot{v}&=&-v^{2}-uv^{2}+u(2+u).
\end{eqnarray}
Figure 11b represents only one orbit of the unit-time map of ($\ref{nilsedlobes1}$) on the $v$-axis.
\begin{center}
\includegraphics[width=4.5cm,bb=0in 0in 3.5in 3.5in]{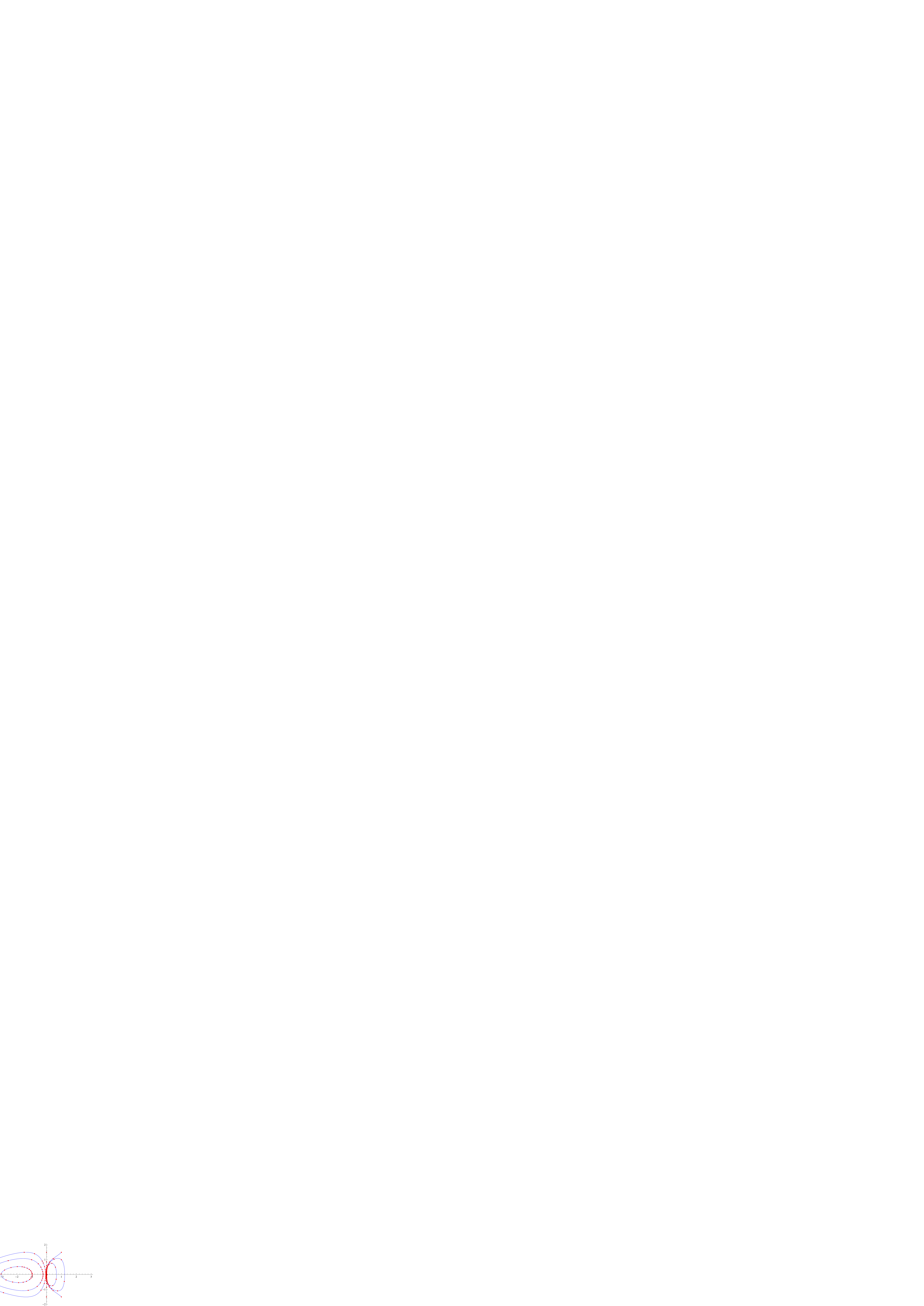} \hskip 1cm  \includegraphics[width=4.5cm,bb=0in 0in 3.5in 3.5in]{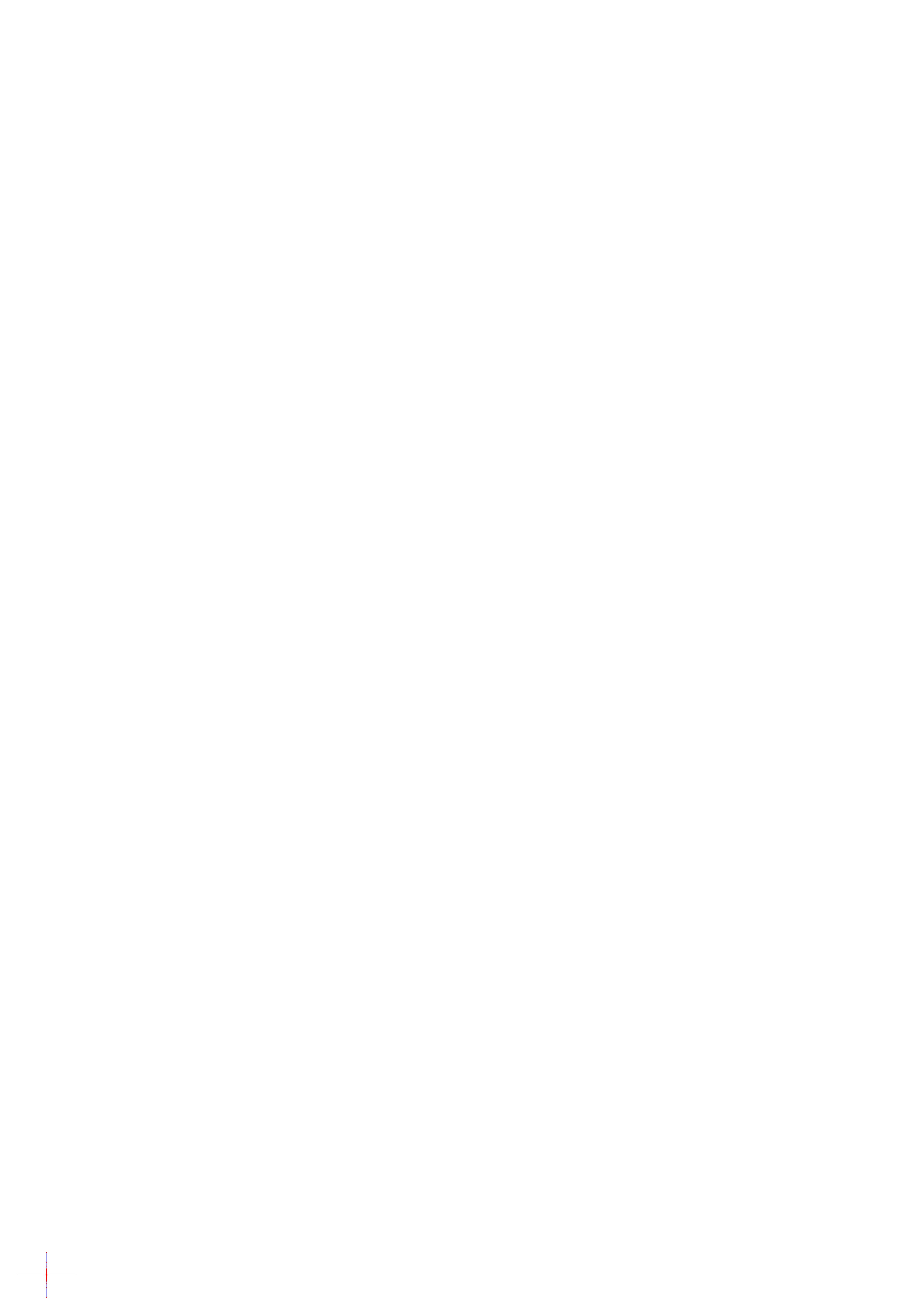}\\ 
\textbf{Figure 11a} system ($\ref{nilsedlobes1}$)  \hskip 1cm \textbf{Figure 11b} $\dim_{B}S(0,v_1)=1/2$
\end{center}

Figures 12a represents the discrete orbits of the unit-time map of the system ($\ref{nilsedlobes}$) with $k=2$
\begin{eqnarray} \label{nilsedlobes2}
\dot{u}&=&-2uv^{3}+u^2v^{3}\nonumber\\
\dot{v}&=&-v^{4}-uv^{4}-u^{2}(2+u)^{2}.
\end{eqnarray}

At Figure 12b we can see only one orbit of the unit-time map of ($\ref{nilsedlobes2}$) on the $v$-axis.

\begin{center}
\includegraphics[width=4.5cm,bb=0in 0in 3.5in 3.5in]{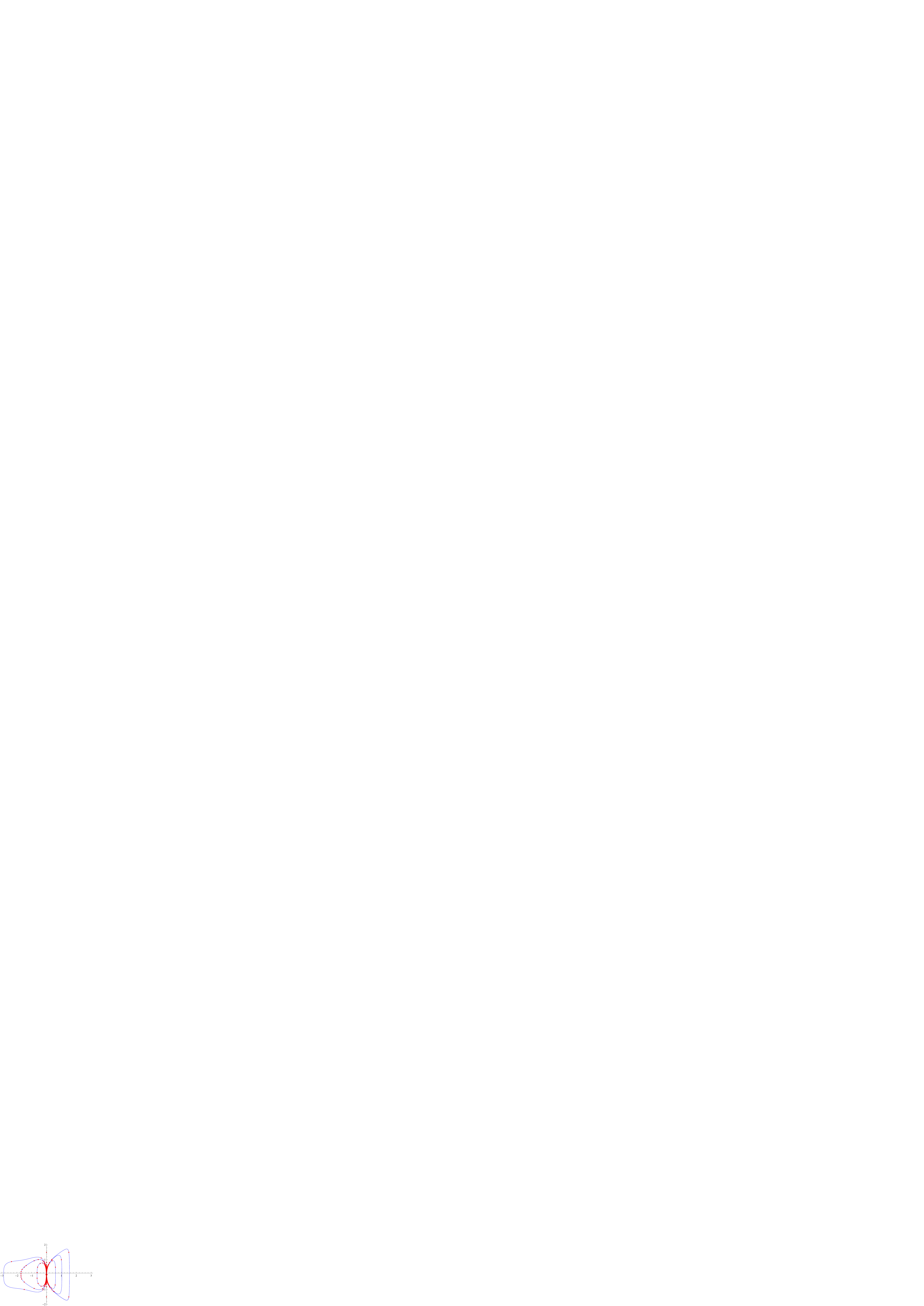} \hskip 1cm  \includegraphics[width=4.5cm,bb=0in 0in 3.5in 3.5in]{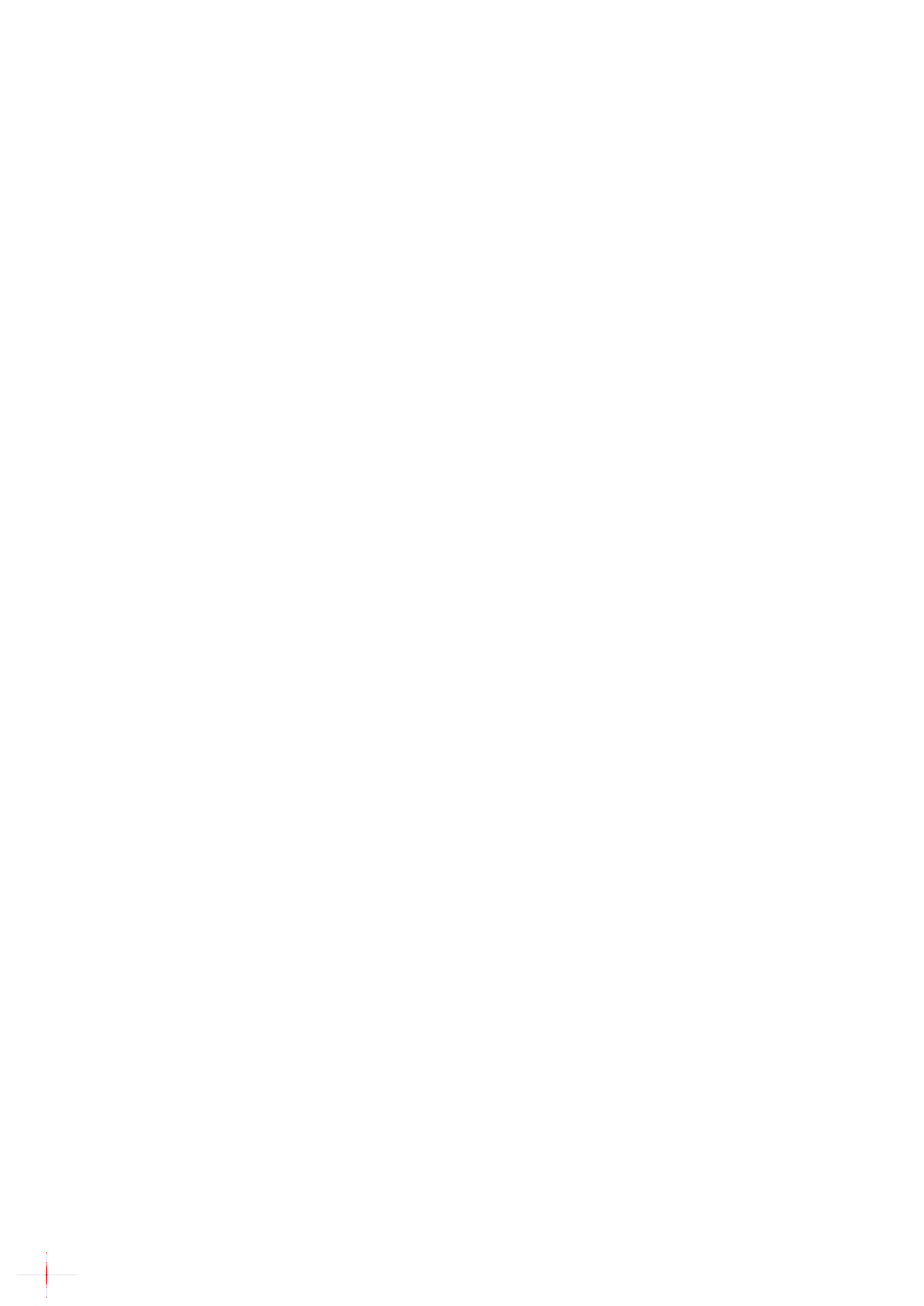}\\ 
\textbf{Figure 12a} system ($\ref{nilsedlobes2}$) \hskip 1cm \textbf{Figure 12b} $\dim_{B}S(0,v_1)=3/4.$
\end{center}


\begin{thebibliography}{99}


\bibitem{aga} M.J.\ Alvarez, A.\  Gasull, {Monodromy and stability for nilpotent critical points}, Int. J. Bif. Chaos, Vol. 15, No. 4 (2005), 1253-1265.

\bibitem{calito} M.\  Caubergh, J.\   Llibre, J.\   Torregrosa, {Global classification of a class of cubic vector
fields whose canonical regions are period annuli}, Int. J. Bif. Chaos, 7 (2011), 1831-1867.


\bibitem{dla} F.\ Dumortier, J.\  Llibre, J.\ C.\  Art$\acute{e}$s, \textit{Qualitative Theory of Planar Differential Systems} (2006), Springer-Verlag Berlin.

\bibitem{neveda} N.\  Elezovi\'c, V.\  \v Zupanovi\'c, D.\  \v Zubrini\'c, Box dimension of trajectories 
of some discrete dynamical systems, Chaos, Solitons \& Fractals Vol.\   34,  2 (2007), 244-252.   

\bibitem{fa} K.\ Falconer, \textit{Fractal Geometry: Mathematical Foundations and Applications}, Chichester: John Wiley and Sons (1990), Chichester.

\bibitem{gi} J.\ Gin\' e, Analytic integrability of nilpotent cubic systems with degenerate infinity, Int. J. Bif. Chaos, Vol. 11, No. 8 (2001), 2299-2304.

 \bibitem{ha} K.\ P.\ Harikrishnan,  R.\ Misra, G.\ Ambika, Revisiting the box counting algorithm for the correlation dimension analysis of hyperchaotic time series, Commun. Nonlinear Sci. Numer. Simulat. 17 (2012) 263-276.

\bibitem{laho} L.\ Horvat Dmitrovi\' c, Box dimension and bifurcations of one-dimensional discrete dynamical systems, Discrete Contin.\   Dyn.\   Syst.\   32 (2012), no.\   4, 1287-1307. 

\bibitem{laho2} L.\ Horvat Dmitrovi\' c, Fractal analysis of Neimark-Sacker bifurcation, preprint (2011).

\bibitem{laho3} L.\ Horvat Dmitrovi\' c, Box dimension of (non)hyperbolic fixed point and singularity of dynamical systems in $\mathbb{R}^{n}$, preprint (2011).

 \bibitem{joyal} P.\ Joyal, Saddle Quantities and Applications, J.\ Differ.\  Eqn.\ 78 (1989), 375-399.
  
\bibitem{kuz} Y.\ A.\ Kuznetsov, \textit{Elements of Applied Bifurcation Theory}, Springer-Verlag New York, USA, (1998).
        
\bibitem{lapo} M.\ L.\ Lapidus, C.\  Pomerance, The Riemann Zeta-function and the one-dimensional Weyl-Berry Conjecture for fractal drums, Proc.\   London Math.\ Soc.\ (3) 66 (1993), no.\ 1, 41-69.


\bibitem{li} W.\  Li, H.\  Wu, Isochronous properties in fractal analysis of some planar vector fields, Bull.\   Sci.\   math.\   134 (2010), 857-873.
\bibitem{mardesic} P.\ Marde\v si\' c, {\em Chebyshev systems and the versal unfolding of the cusp of order $n$} (1998), Hermann, \' Editeurs des Sciences et des Arts, Paris.
\bibitem{mrz} P.\   Marde\v si\' c, M.\   Resman, V.\  \v Zupanovi\' c,  Multiplicity of fixed points and $\varepsilon-$neighborhoods of orbits, arXiv:1108.4707 (2011).

\bibitem{pa} M.\ Pa\v si\' c, Minkowski-Bouligand dimension of solutions of the one-dimensional $p$-Laplacian,
        J.\ Differ.\ Eqn.\ 190 (2003), 268-305.
        
\bibitem{pa1} M.\ Pa\v si\' c, Fractal oscillations for a class of second-order linear differential equations of Euler type, J.\   Math.\   Anal.\   Appl.\   341 (2008), 211-223.

\bibitem{pazuzu} M.\ Pa\v si\' c, D.\   \v Zubrini\' c, V.\  \v Zupanovi\' c, Oscillatory and phase dimensions of solutions of some second-order differential equations, Bull.\   sci.\   math.\   133 (8) (2009), 859-874.

\bibitem{bogdanov} S.\ P\'erez Gonz\'ales, A.\ Gasull, J.\ Torregrosa, Global study of the Bogdanov-Takens
bifurcation curve, preprint (2011)

\bibitem{rzz} G.\ Radunovi\'c, D.\ \v Zubrini\'c, and V.\  \v Zupanovi\'c, Fractal analysis of Hopf bifurcation at infinity, to appear in Int. J. Bif. Chaos (2012).
\bibitem{r} M.\ Resman, Formal classification of parabolic diffeomorphisms and asymptotic development of $\varepsilon$-neighborhoods of orbits, preprint (2012).
\bibitem{rezz} M.\  Resman, D.\   \v Zubrini\'c, V.\  \v Zupanovi\'c,  Poincar\'e map of a class of degenerate foci and applications, preprint (2011).

\bibitem{sz} E.\ Str\' o\.zyna, H.\ \.Zo\l adek, The analytic and formal normal form for the nilpotent singularity, J.\ Differ.\ Eqn.\ 179, (2002), 479-537.

\bibitem{takens} F.\ Takens, Unfoldings of certain singularities of vector fields: Generalized Hopf bifurcations, J.\ Differ.\  Eqn.\  14 (1973) 476-493.
\bibitem{t} C.\   Tricot, \textit{Curves and Fractal Dimension} (1995), Springer-Verlag New York. 
\bibitem{zuzu} D.\  \v Zubrini\'c, V.\  \v Zupanovi\'c, Fractal analysis of spiral trajectories of some
planar vector fields, Bull.\   sci.\   math.\   129/6 (2005), 457-485.


 \bibitem{wu} Y.\ Wu, P.\ Li, H.\ Chen, Center conditions and bifurcation of limit cycles at three-order
nilpotent critical point in a cubic Lyapunov system, Commun. Nonlinear Sci. Numer. Simulat. 17 (2012), 292-304.

\bibitem{zuzu3} D.\  \v Zubrini\'c, V.\  \v Zupanovi\'c,
Fractal analysis of spiral trajectories of some vector fields in ${\mathbb R}^3$, 
C.\   R.\   Acad.\   Sci.\   Paris, S\'erie I, Vol.\   342, 12 (2006), 959-963.

\bibitem{zuzu3q} D.\ \v Zubrini\'c, V.\ \v Zupanovi\'c, Box dimension of spiral trajectories of some
vector fields in $\eR^3$, Qual.\ Theory Dyn.\ Syst.\ Vol 6 (2005), 251-272.

\bibitem{belg} D.\  \v Zubrini\'c, V.\  \v Zupanovi\'c, Poincar\'e map in fractal analysis of 
spiral trajectories of planar vector fields, Bull.\   Belg.\   Math.\   Soc.\   Simon Stevin, 15 (2008) 947-960.
        
    \bibitem{z} V.\  \v Zupanovi\' c, Topological Equivalence of Planar Vector Fields and Their Generalised Principal Part, J.\  Differ.\  Eqn.\ 167, (2000), 1-15.

\bibitem{zuzu4} V.\  \v Zupanovi\' c, D.\  \v Zubrini\' c, \textit{Fractal dimension in dynamics}, Encyclopedia of Math.\   Physics, J.-P.\   Fran\c coise, G.\ L.\   Naber, S.\ T.\   Tsou (Eds.\ ), vol.\ 2 (2006) Elsevier, Oxford. 
    
 
\end{thebibliography}
\end{document}